\documentclass[11pt]{amsart}
\usepackage[all]{xy}
\usepackage{mathrsfs}
\usepackage{amssymb,amsmath,amsgen,amsxtra,amsfonts} 
\usepackage{dsfont,times}

\begin{document}
%
%   D e f i n i t i o n s
%
% gibt's die folgenden Definition schon?
%\newcommand{\qed}{\mbox{}\hfill$\Box$}
%\newcommand{\mod}{{\:\mbox{mod}\:}}
%
\theoremstyle{definition}
\newtheorem{Definition}{Definition}[section]
\newtheorem*{Definitionx}{Definition}
\newtheorem{Convention}{Definition}[section]
\newtheorem{Construction}{Construction}[section]
\newtheorem{Example}[Definition]{Example}
\newtheorem{Examples}[Definition]{Examples}
\newtheorem{Remark}[Definition]{Remark}
\newtheorem*{Remarkx}{Remark}
\newtheorem{Remarks}[Definition]{Remarks}
\newtheorem{Caution}[Definition]{Caution}
\newtheorem{Conjecture}[Definition]{Conjecture}
\newtheorem*{Conjecturex}{Conjecture}
\newtheorem{Question}[Definition]{Question}
\newtheorem{Questions}[Definition]{Questions}
\newtheorem*{Acknowledgements}{Acknowledgements}
\newtheorem*{Disclaimer}{Disclaimer}
\theoremstyle{plain}
\newtheorem{Theorem}[Definition]{Theorem}
\newtheorem*{Theoremx}{Theorem}
\newtheorem{Proposition}[Definition]{Proposition}
\newtheorem*{Propositionx}{Proposition}
\newtheorem{Lemma}[Definition]{Lemma}
\newtheorem{Corollary}[Definition]{Corollary}
\newtheorem*{Corollaryx}{Corollary}
\newtheorem{Fact}[Definition]{Fact}
\newtheorem{Facts}[Definition]{Facts}
\newtheoremstyle{voiditstyle}{3pt}{3pt}{\itshape}{\parindent}%
{\bfseries}{.}{ }{\thmnote{#3}}%
\theoremstyle{voiditstyle}
\newtheorem*{VoidItalic}{}
\newtheoremstyle{voidromstyle}{3pt}{3pt}{\rm}{\parindent}%
{\bfseries}{.}{ }{\thmnote{#3}}%
\theoremstyle{voidromstyle}
\newtheorem*{VoidRoman}{}

% abgeschrieben aus The LaTeX Companion, 2nd edition,
% von Mittelback & Goossens
%
\newcommand{\prf}{\par\noindent{\sc Proof.}\quad}
\newcommand{\blowup}{\rule[-3mm]{0mm}{0mm}}
\newcommand{\cal}{\mathcal}
\newcommand{\Aff}{{\mathds{A}}}
\newcommand{\BB}{{\mathds{B}}}
\newcommand{\CC}{{\mathds{C}}}
\newcommand{\EE}{{\mathds{E}}}
\newcommand{\FF}{{\mathds{F}}}
\newcommand{\GG}{{\mathds{G}}}
\newcommand{\HH}{{\mathds{H}}}
\newcommand{\NN}{{\mathds{N}}}
\newcommand{\ZZ}{{\mathds{Z}}}
\newcommand{\PP}{{\mathds{P}}}
\newcommand{\QQ}{{\mathds{Q}}}
\newcommand{\RR}{{\mathds{R}}}
\newcommand{\Liea}{{\mathfrak a}}
\newcommand{\Lieb}{{\mathfrak b}}
\newcommand{\Lieg}{{\mathfrak g}}
\newcommand{\Liem}{{\mathfrak m}}
\newcommand{\ideala}{{\mathfrak a}}
\newcommand{\idealb}{{\mathfrak b}}
\newcommand{\idealg}{{\mathfrak g}}
\newcommand{\idealm}{{\mathfrak m}}
\newcommand{\idealp}{{\mathfrak p}}
\newcommand{\idealq}{{\mathfrak q}}
\newcommand{\idealI}{{\cal I}}
\newcommand{\lin}{\sim}
\newcommand{\num}{\equiv}
\newcommand{\dual}{\ast}
\newcommand{\iso}{\cong}
\newcommand{\homeo}{\approx}
\newcommand{\mm}{{\mathfrak m}}
\newcommand{\pp}{{\mathfrak p}}
\newcommand{\qq}{{\mathfrak q}}
\newcommand{\rr}{{\mathfrak r}}
\newcommand{\pP}{{\mathfrak P}}
\newcommand{\qQ}{{\mathfrak Q}}
\newcommand{\rR}{{\mathfrak R}}
%
%  evtl. auch \"uber \mathbb oder \Bbb
%
\newcommand{\dq}{{``}}
\newcommand{\OO}{{\cal O}}
\newcommand{\numero}{{n$^{\rm o}\:$}}
\newcommand{\into}{{\hookrightarrow}}
\newcommand{\onto}{{\twoheadrightarrow}}
\newcommand{\Spec}{{\rm Spec}\:}
\newcommand{\BigSpec}{{\rm\bf Spec}\:}
\newcommand{\Spf}{{\rm Spf}\:}
\newcommand{\Proj}{{\rm Proj}\:}
\newcommand{\Pic}{{\rm Pic }}
\newcommand{\Br}{{\rm Br}}
\newcommand{\NS}{{\rm NS}}
\newcommand{\Sym}{{\mathfrak S}}
\newcommand{\Aut}{{\rm Aut}}
\newcommand{\Autp}{{\rm Aut}^p}
\newcommand{\Hom}{{\rm Hom}}
\newcommand{\ord}{{\rm ord}}
\newcommand{\divisor}{{\rm div}}
\newcommand{\Brauer}{\widehat{{\rm Br}}}
\newcommand{\ModulPolAmple}{{\mathscr M}_{\rm CV, ample}}
\newcommand{\ModulPolSmooth}{{\mathscr M}_{\rm CV, smooth}}
\newcommand{\Modul}{{\mathscr M}_{\rm Enriques}}
\newcommand{\MMod}{{\mathscr M}}
\newcommand{\ModulGrp}{{\mathscr M}_{\rm Grp, 2}}
\newcommand{\Def}{{\rm Def}}
\newcommand{\Defp}{{\rm Def}^p}
\newcommand{\DefGrp}{{\rm Def}_{{\rm Grp}, 2}}
\newcommand{\piet}{{\pi_1^{\rm \acute{e}t}}}
\newcommand{\Het}[1]{{H_{\rm \acute{e}t}^{{#1}}}}
\newcommand{\Hcris}[1]{{H_{\rm cris}^{{#1}}}}
\newcommand{\HdR}[1]{{H_{\rm dR}^{{#1}}}}
\newcommand{\hdR}[1]{{h_{\rm dR}^{{#1}}}}
\newcommand{\defin}[1]{{\bf #1}}

\title[Enriques Surfaces]{Arithmetic Moduli and Lifting of Enriques Surfaces}

\subjclass[2000]{14J28, 14J10}
%\keywords{Enriques surface, K3 surface}

\author{Christian Liedtke}
\address{TU M\"unchen, Zentrum Mathematik - M11, Boltzmannstr. 3, D-85748 Garching bei M\"unchen, Germany}
\curraddr{}
\email{liedtke@ma.tum.de}

\date{July 13, 2013}

\dedicatory{in memoriam Torsten Ekedahl}
\maketitle

\begin{abstract}
  We construct the moduli space of Enriques surfaces in positive characteristic and
  eventually over the integers, 
  and determine its local and global structure.
  As an application, we show lifting of Enriques surfaces to characteristic zero.
  The key observation is that the canonical double
  cover of an Enriques surface is birational to the 
  complete intersection of three quadrics in $\PP^5$, even
  in characteristic $2$.
\end{abstract}

\section*{Introduction}

In order to give examples of algebraic surfaces with 
$h^1(\OO_X)=h^2(\OO_X)=0$ that are not rational, 
Enriques constructed the first Enriques 
surfaces at the end of the 19th century.
From the point of view of the Kodaira--Enriques classification, 
these surfaces form one of the four classes
of minimal surfaces of Kodaira dimension zero.
More precisely, these four classes consist of Abelian surfaces,
K3 surfaces, Enriques surfaces and 
(Quasi-)Hyperelliptic surfaces.

In characteristic $\neq2$, Enriques surfaces behave extremely nice:
deformations are unobstructed by results of Illusie \cite{Illusie}
and Lang \cite{Lang}.
Over the complex numbers, their moduli space is irreducible, 
unirational and $10$-dimensional, and
Kond\={o} \cite{Kondo} showed even rationality. 
Next, their fundamental groups are of order $2$ and their
universal covers are K3 surfaces.
Moreover, Cossec \cite{Cossec Models} and Verra \cite{Verra}
found explicit equations of these K3-covers:
they are birational to complete intersections 
of three quadrics in $\PP^5$ and the $\ZZ/2\ZZ$-action
can be written down explicitly.
For generic Enriques surfaces this was already known
to Enriques himself \cite{Enriques}.
Finally, 
Enriques surfaces in characteristic $\neq2$ lift over 
the Witt ring, which is due to Lang \cite{Lang}.

In characteristic $2$, the situation is more complicated:
first of all, as shown by Bombieri and Mumford \cite{bm3},
Enriques surfaces fall into three classes, called
{\em classical}, {\em singular} and {\em supersingular}.
Although they still possess canonically defined flat double covers, 
which have trivial dualizing sheaves and which ``look'' cohomologically
like K3 surfaces,
these are in general only integral Gorenstein surfaces and
may not be normal.
It also happens that deformations are
obstructed and finally, supersingular Enriques surfaces do not
lift over the Witt ring.
\medskip

In this article, we clarify the situation in positive characteristic,
and especially in characteristic $2$.
We start with the description of the canonical double cover.

\begin{Theoremx}
  Let $\pi:\widetilde{X}\to X$ be the K3(-like) canonical double
  cover of an Enriques surface $X$.
  Then, there exists a morphism $\varphi$
  $$
     \widetilde{X}\,\to\,\varphi(\widetilde{X})\,\subseteq\,\PP^5
  $$
  that is birational onto its image. 
  The image $\varphi(\widetilde{X})$ is a complete intersection of
  three quadrics.
\end{Theoremx}

The exceptional locus of $\varphi$ is, in a certain sense,
a union of $ADE$-curves, see Theorem \ref{three quadrics} for a precise
statement.
Next, $\pi$ is a torsor under the 
finite flat group scheme $G:=\Pic^\tau(X)^D$, which is of length 
$2$ (we refer to Section \ref{sec:generalities} for definition).
We describe the linear $G$-action on $\PP^5$ induced by $\varphi$, 
as well as the equations of the $G$-invariant quadrics cutting out 
$\varphi(\widetilde{X})$ in Proposition \ref{three quadrics explicit}.
As a consequence of this structure result, 
we conclude the following (see also Remark \ref{reid bombieri mumford}):

\begin{Corollaryx}
 All Enriques surfaces in arbitrary characteristic arise
 via the Bombieri--Mumford--Reid construction in \cite[\S3]{bm3}.
\end{Corollaryx}

Next, we study polarized families of Enriques surfaces,
which turn out to be interesting in their own right.
In a sense made precise below, the moduli space of such polarized surfaces
behaves as one would hope the moduli space to look like, whereas the
moduli space of unpolarized Enriques surfaces displays some pathologies
(see below or Remark \ref{resolution bad}).

As polarizations, we choose the following class of invertible sheaves,
which is rather natural and minimal 
(see the beginning of Section \ref{sec:three quadrics}):

\begin{Definitionx}
 A {\em Cossec--Verra polarization} on an Enriques surface $X$ is an
 invertible sheaf ${\cal L}$ with ${\cal L}^2=4$  and
 such that every genus-one fibration $|2E|$ on $X$  satisfies
 $\deg {\cal L}|_E\geq2$.
\end{Definitionx}

In general, Cossec--Verra polarizations need not be ample, but are merely
big and nef. 
Thus, to call them quasi-polarizations might be more appropriate, but for sake
of readability we have decided not to.
On the other hand, every Enriques surface 
carries such a polarization.
This is different from algebraic K3 surfaces, all of
which are polarizable, but where we need infinitely many
types of polarizations to capture every one of them.

In general, Cossec--Verra polarizations are not unique, and we
refer to Proposition \ref{prop:polarization} for 
quantitative results.
Contracting those curves that have zero-intersection
with a given Cossec--Verra polarization $\cal L$, we obtain a pair $(X',{\cal L}')$, where
$X'$ is an Enriques surface with at worst Du~Val singularities
and ${\cal L}'$ is an {\em ample} Cossec--Verra polarization.

\begin{Theoremx}
 Let $X'$ be an Enriques surface over $k$
 with at worst Du~Val singularities
 admitting an ample Cossec--Verra polarization.
 \begin{enumerate}
  \item If $X'$ is not supersingular then it lifts over the Witt ring $W(k)$.
  \item If $X'$ is supersingular then it lifts over $W(k)[\sqrt{2}]$,
    but not over $W(k)$. 
 \end{enumerate}
\end{Theoremx}

Using Artin's results \cite{Artin Brieskorn} on simultaneous resolutions 
of singularities in families, we conclude that for every Enriques
surface $X$ over $k$, there exists an algebraic space that is smooth
over a possibly ramified extension of $W(k)$ with special
fiber $X$.
Put a little bit sloppily, we conclude:

\begin{Corollaryx}
  Enriques surfaces lift to characteristic zero.
\end{Corollaryx}

We refer to Theorem \ref{thm:general lifting} for precise results, and note that
lifting of Enriques surfaces in characteristic $\neq2$ was already established by Lang \cite{Lang}.

Next, we construct and study the moduli space $\ModulPolAmple$, whose geometric points
correspond to pairs $(X,{\cal L})$ where $X$ is an Enriques surface with at 
worst Du~Val-singularities,
and $\cal L$ is an ample Cossec--Verra polarization.
This moduli space behaves extremely nice, even in characteristic $2$:

\begin{Theoremx}
  $\ModulPolAmple$ is a quasi-separated Artin stack of finite type over $\ZZ$.
  \begin{enumerate}
   \item If $p\neq2$ then $\ModulPolAmple\otimes_\ZZ\FF_p$ is irreducible, unirational, $10$-dimensional, 
     smooth over $\FF_p$, and even a Deligne--Mumford stack.
   \item On the other hand, $\ModulPolAmple\otimes_\ZZ\FF_2$ consists of two components
    $$
       \MMod^{\mu_2} \mbox{ \quad and \quad }\MMod^{\ZZ/2\ZZ}\,.
     $$  
     both of which are irreducible, unirational, smooth, and $10$-dimensional Artin stacks over $\FF_2$.
     More precisely,
     \begin{enumerate}
     \item[-] they intersect transversally along an irreducible, unirational, smooth and $9$-dimensional 
        closed substack $\MMod^{\alpha_2}$, where
     \item[-] $\MMod^{\alpha_2}$ parametrizes supersingular Enriques surfaces, and
     \item[-] $\MMod^G\,-\,\MMod^{\alpha_2}$ parametrizes
           singular Enriques surfaces ($G=\mu_2$) and 
           classical Enriques surfaces ($G=\ZZ/2\ZZ$), respectively.
     \end{enumerate}
  \end{enumerate}
\end{Theoremx}

The lifting result and the description of $\ModulPolAmple$ give a beautiful
picture, and one might expect a similar picture for moduli spaces of smooth Enriques surfaces, polarized as well
as unpolarized. Unfortunately, this is not the case:
First of all, the stack $\Modul$ of smooth and unpolarized Enriques surfaces is not quasi-separated (that is,
its diagonal is not quasi-compact), making geometric statements delicate, see
Remark \ref{rem:not quasi-separated}.
But even putting aside this technical issue, deformations of classical Enriques surfaces in characteristic $2$
may be obstructed, and deformations of supersingular Enriques in characteristic $2$ may be more obstructed
than the normal crossing situation along $\MMod^{\alpha_2}$ of the above theorem.
This was observed by Ekedahl and Shepherd-Barron in \cite{Ekedahl} and led them to
introduce {\em exceptional Enriques surfaces}.
%In \cite{Ekedahl unpublished}, the hulls of the deformation functors of such surfaces were computed.

To explain these obstructions, we consider the moduli space $\ModulPolSmooth$ of smooth Enriques
surfaces with Cossec--Verra polarizations, which is in fact a quasi-separated Artin stack of finite
type over $\ZZ$.
Then, we obtain a diagram 
$$
 \xymatrixcolsep{3pc}
 \xymatrix{
      & \ar[dl]_{\Phi_{\rm cont}} \ar[dr]^{\Phi_{\rm forget}} \ModulPolSmooth\\
      \ModulPolAmple & & \Modul }
$$
where $\Phi_{\rm forget}$ is the functor that forgets the Cossec--Verra polarization, 
and where $\Phi_{\rm cont}$ is the functor that contracts curves that have zero-intersection with
the Cossec--Verra polarization (see Section \ref{sec:final} for precise definitions).
Now, although $\Phi_{\rm cont}$ is a bijection on geometric points,
it is here that $\ModulPolSmooth$ picks up singularities, and we refer to  
Remark \ref{resolution bad} for details.
We also refer the reader to Section \ref{sec:final} for structure results on
$\Modul$ and $\ModulPolSmooth$.
\medskip

This article is organized as follows:

After reviewing a couple of general facts in Section \ref{sec:generalities},
we study projective and birational models of the canonical double 
cover of Enriques surfaces in Section \ref{sec:projective models}.
This extends work of Cossec \cite{Cossec Models} to
characteristic $2$.
The main difficulty is that Saint-Donat's results \cite{Saint-Donat} on
linear systems on K3 surfaces cannot be applied to this double cover and
we have to find rather painful ways around.

In Section \ref{sec:three quadrics} we show that the canonical double
cover of an Enriques surface is birational to the complete intersection 
of three quadrics in $\PP^5$.
We introduce the notion of a Cossec--Verra polarization and establish
a couple of general facts about them.
Finally, we explicitly describe the action of the finite flat group scheme 
of length two, which acts on this complete intersection.
In particular, we will see that every Enriques surface arises via the
Bombieri--Mumford--Reid construction of \cite{bm3}.

In Section \ref{sec:moduli} we study pairs of Enriques surfaces
with at worst Du~Val singularities together with ample Cossec--Verra
polarizations.
We prove their lifting to characteristic zero
and construct the moduli space $\ModulPolAmple$ of such pairs.
Using the results of Section \ref{sec:three quadrics}, all
boils down to describing deformations of complete 
intersections of three quadrics in $\PP^5$
together with the action of a finite flat group scheme of length $2$.

Finally, in Section \ref{sec:final}, we relate $\ModulPolAmple$
to moduli spaces of smooth Enriques surfaces, both polarized and
unpolarized.
These are connected via Artin's functor of simultaneous resolutions
of singularities and it is here, that these moduli spaces pick up
singularities.

\begin{Acknowledgements}
 It is pleasure for me  to thank Brian~Conrad, Igor~Dolgachev, David~Eisenbud,  
 Jack~Hall, Sven~Meinhardt, and especially Torsten~Ekedahl 
 for discussions and help.
 Also, I thank the referee for remarks and comments.
 Moreover, I gratefully acknowledge funding from DFG under 
 research grant LI 1906/1-1 and Transregio SFB 45 and thank the departments
 of mathematics at Stanford university and Bonn university for hospitality, 
 where parts of this article were written.
\end{Acknowledgements}

\section{Generalities on Enriques Surfaces}
\label{sec:generalities}

We start by recalling a couple of general facts on 
Enriques surfaces, and refer to 
\cite[\S3]{bm3} and \cite[Chapter I]{Cossec; Dolgachev}
for details, proofs, and further references.

\begin{Definition}
  A (smooth) {\em Enriques surface} is a smooth and proper surface
  $X$ of finite type over an algebraically closed field
  $k$ of characteristic $p\geq0$ such that
  $$
  \omega_X\,\equiv\,\OO_X \mbox{ \quad and \quad }b_2(X)\,=\,10,
  $$ 
  where $\equiv$ denotes numerical equivalence, and
  $b_i$ denotes the $i$.th \'etale or crystalline Betti number.
\end{Definition}

If $X$ is an Enriques surface, then we have 
$$
\chi(\OO_X)\,=\,1\mbox{ \quad and \quad }b_1(X)=0\,.
%\mbox{ \quad and \quad }{\rm rank}\,\Pic(X)\,=\,b_2(X)\,. 
$$
Moreover, $\omega_X^{\otimes2}\iso\OO_X$ holds in every characteristic,
and $\omega_X\not\iso\OO_X$ if and only if $h^1(\OO_X)=0$.
In characteristic $p\neq2$ we have $h^1(\OO_X)=0$, 
whereas for $p=2$ only the inequality $h^1(\OO_X)\leq1$ holds true. 
Thus, if $H^1(\OO_X)$ is non-zero, it makes sense to study the 
action of the absolute Frobenius $F$ on it, which is
either zero or a bijection.
Combing these results, we have the following definition and characterization, 
due to Bombieri and Mumford \cite[\S3]{bm3}.

\begin{Definition}
 \label{def: enriques types}
 An Enriques surface $X$ is called 
 \begin{enumerate}
  \item {\em classical} if $h^1(\OO_X)=0$, hence $\omega_X\not\iso\OO_X$
   and $\omega_X^{\otimes 2}\iso\OO_X$,
  \item {\em singular} if $h^1(\OO_X)=1$, hence $\omega_X\iso\OO_X$ and
   $F$ is bijective on $H^1(\OO_X)$,
  \item {\em supersingular} if $h^1(\OO_X)=1$, hence $\omega_X\iso\OO_X$ and
   $F$ is zero on $H^1(\OO_X)$.
 \end{enumerate}
\end{Definition}

The Picard scheme of $X$ is smooth if and only if it is classical.
Moreover, $\Pic^\tau(X)$ is a group scheme of length $2$.
Let us recall, for example from \cite{Tate; Oort}, that a finite group scheme
of length $2$ over an algebraically closed field of characteristic $p\neq2$ is
isomorphic to $\ZZ/2\ZZ$, whereas there are three
such group schemes for $p=2$, namely $\ZZ/2\ZZ$, $\mu_2$,
and $\alpha_2$.
In fact, these three group schemes correspond to
the three classes of Enriques surfaces in Definition \ref{def: enriques types}:
$\Pic^\tau(X)$ is isomorphic to  $\ZZ/2\ZZ$ if $X$ is classical, 
to $\mu_2$ if $X$ is singular, and to $\alpha_2$ if $X$ is supersingular.

In any case,
$\Pic^\tau$ gives rise to finite and flat morphism of degree $2$
$$
   \pi\,:\,\widetilde{X}\,\longrightarrow\,X,
$$
which is a torsor under $G:=(\Pic^\tau(X))^D$, and
where $-^D={{\cal H}om}(-,\GG_m)$ denotes
Cartier duality.
%We refer to  \cite[Proposition (6.2.1)]{Raynaud} for details.

In particular, if $p\neq2$ or if $X$ is a singular Enriques surface, then
$G\iso\ZZ/2\ZZ$, the morphism $\pi$ is an \'etale Galois cover of degree $2$
and $\widetilde{X}$ is a K3 surface.
In the remaining cases, that is, $X$ is classical with $p=2$ or supersingular, 
$\pi$ is purely inseparable and $\widetilde{X}$
is never smooth, possibly even non-normal.
In any case, $\widetilde{X}$ is an integral Gorenstein surface 
with invariants
$$
\omega_{\widetilde{X}}\,\iso\,\OO_{\widetilde{X}},\mbox{ \quad }
\chi(\OO_{\widetilde{X}})\,=\,2\mbox{ \quad and \quad }
h^1(\OO_{\widetilde{X}})\,=\,0,
$$
that is, ``K3-like''.
Having only an integral Gorenstein surface rather than a smooth 
K3 surface as double cover, is one of the main reasons why Enriques surfaces
in characteristic $2$ are so difficult to come by.

Finally, we recall some numerical invariants, and refer to
\cite[Theorem 0.11]{Lang}, \cite[Section II.7.3]{Illusie}, and
\cite[Proposition 1.4.2]{Cossec; Dolgachev} for details.
In the following table, $\Theta_X$ denotes the tangent sheaf.
$$
\begin{array}{cl|cc|ccc}
p&\mbox{ type } &h^1(\OO_X)&h^0(\Omega_X^1)&h^0(\Theta_X)&h^1(\Theta_X)&h^2(\Theta_X)\\
\hline
2 & \mbox{ classical }     & 0 & 1 & a & 10+2a & a\\
  & \mbox{ singular }      & 1 & 0 & 0 & 10 & 0 \\
  & \mbox{ supersingular } & 1 & 1 & 1 & 12 & 1 \\
\hline
\neq2 &                    & 0 & 0 & 0 & 10 & 0
\end{array}
$$
Classical Enriques surfaces satisfy $a\leq1$ and surfaces with
$a=1$ have been described and explicitly classified by
Ekedahl and Shepherd-Barron \cite{Ekedahl} and
Salomonsson \cite{Salomonsson}.
The existence of classical Enriques surfaces with $a=1$ is quite an
unpleasant surprise,
to which we come back in Remark \ref{resolution bad}.

\section{Projective models of the canonical double cover}
\label{sec:projective models}

In this section we study linear systems and 
projective models of the canonical double cover $\widetilde{X}$ 
of an Enriques surface $X$.
Since there is no canonical polarization on $\widetilde{X}$, 
the best thing 
to do is to consider pull-backs of invertible sheaves from $X$
with positive self-intersection number. 
In characteristic $\neq2$, this has been carried out by
Cossec 
\cite[Section 8]{Cossec Models}: % and \cite[Chapter IV]{Cossec; Dolgachev}.
such a pull-back defines a morphism that is
either birational onto its image 
or generically finite of degree $2$ onto a rational surface.
%Which of these cases occurs is controlled by the numerical
%invariant $\Phi$ (see below).
Cossec's proof relies on Saint-Donat's 
analysis \cite{Saint-Donat} of linear systems on K3 surfaces
that he applies to $\widetilde{X}$.
In characteristic $2$, the main difficulty is that $\widetilde{X}$
is in general only an integral Gorenstein surface, and so we have
to take rather painful detours.

As before, we denote by $\pi:\widetilde{X}\to X$ the canonical double cover
of an Enriques surface $X$.
For an effective divisor $C$ on $X$, we pull back $\OO_X(C)$ 
to $\widetilde{X}$ and study the associated (rational) map
from $\widetilde{X}$ to projective space.
To do so, we consider the ``positivity measure'' $\Phi$, 
% for $\OO_X(C)$ 
introduced in \cite[Chapter III \S2]{Cossec; Dolgachev}:

\begin{Definition}
 For an effective divisor $C$ on an Enriques surface $X$, we set
 $$
   \Phi(C)\,:=\,\frac{1}{2}\,\inf\big\{ E\cdot C,
  \mbox{ where }|E|\mbox{ is a genus one pencil on $X$}\big\} \,.
  $$
\end{Definition}

We note that every Enriques surface possesses at least one genus-one fibration,
and that every divisor defining a genus-one fibration is $2$-divisible
in $\Pic(X)$, % which leads to the notion of ``half-pencils'',
see \cite[Chapter V \S7]{Cossec; Dolgachev}.
In particular, $\Phi(C)$ is a well-defined and non-negative integer.
For example, if $C$ is an irreducible curve with $C^2>0$ then the linear
system $|C|$ is basepoint-free if and only if $\Phi(C)\geq2$, see
\cite[Theorem 4.4.1]{Cossec; Dolgachev}.
As we shall see now and in Theorem \ref{sd3} below,
$\Phi$ also controls the behavior of linear systems
on $\widetilde{X}$.
More precisely:

\begin{Theorem}
 \label{sd1}
% Let $X$ be an Enriques surface and
 Let $C$ be an irreducible curve with $C^2>0$ and $\Phi(C)\geq2$
 on an Enriques surface $X$.
 Then,
 \begin{enumerate}
  \item $C^2\geq4$,
  \item the invertible sheaf $\pi^\ast\OO_X(C)$ on $\widetilde{X}$
    is globally generated,
  \item a generic Cartier divisor in $|\pi^\ast\OO_X(C)|$ is an
    integral Gorenstein curve, which is not hyperelliptic, and
  \item $|\pi^\ast\OO_X(C)|$ gives rise to a morphism 
    $$
   \varphi\,:\,\widetilde{X}\,\longrightarrow\,\PP^{(1+C^2)},
   $$
   which is birational onto an integral surface of degree $2C^2$.
\end{enumerate}
\end{Theorem}

\prf
By \cite[Theorem 4.4.1]{Cossec; Dolgachev},
the linear system $|C|$ is basepoint-free.
In particular, $\pi^\ast\OO_X(C)$ is globally generated on $\widetilde{X}$.
From the formula for $h^0$ of \cite[Corollary 1.5.1]{Cossec; Dolgachev} 
it follows that a generic divisor in $|C|$ is reduced.
Now, if we had $C^2=2$ then 
$|C|$ would define a morphism onto $\PP^1$.
In particular, $C$ would be the class of a fiber, it would satisfy
$C^2=0$, which contradicts $C^2\neq0$, and thus, we conclude
$C^2\geq4$.

Next, we consider the short exact sequence
\begin{equation}
 \label{ses}
 0\,\to\,\OO_X(C)\,\to\,\pi_\ast\pi^\ast\OO_X(C)\,\to\,
  \omega_X(C)\,\to\,0\,.
\end{equation}
We have $h^1(X,\OO_X(C))=0$ by \cite[Theorem 1.5.1]{Cossec; Dolgachev},
which, together with 
\cite[Corollary 1.5.1]{Cossec; Dolgachev} implies
$h^0(\widetilde{X},\pi^\ast\OO_X(C))=2+C^2$.
Thus, $\pi^\ast\OO_X(C)$ gives rise to a morphism 
$\varphi$ from $\widetilde{X}$
to $(1+C^2)$-dimensional projective space.
Also, since the image of $|C|$ is a surface, the
same is true for $\varphi$.
Moreover, $\varphi(\widetilde{X})$
is an integral surface, that is, reduced and irreducible,
since $\widetilde{X}$ is.

If $p\neq2$ or if $X$ is a singular Enriques surface,
then $\widetilde{X}$ is a smooth K3 surface and we compute
$(\pi^\ast C)^2=2C^2$.
Since $\pi^\ast\OO_X(C)$ is globally generated,
we find $2C^2=\deg\varphi\cdot\deg\varphi(\widetilde{X})$.
Since non-degenerate and integral surfaces in $\PP^N$
have degree at least $N-1$ (see \cite[Proposition 0]{eh} for a proof that
is valid in arbitrary characteristic),
we conclude $\deg\varphi\leq2$. 

If $p=2$ and $X$ is classical or supersingular, then
$\pi$ is a torsor under $\mu_2$ or $\alpha_2$.
In particular, $\pi$ is purely inseparable of degree $2$,
and the extension $k(X)\subset k(\widetilde{X})$ of function fields
is obtained by adjoining a square root.
If we denote by $k(X)^{1/2}$
the field that is obtained by adjoining all square roots of $k(X)$,
then the resulting  field extension $k(X)\subset k(X)^{1/2}$ is
purely inseparable.
Moreover, we have an inclusion of fields
$k(X)\subset k(\widetilde{X})\subset k(X)^{1/2}$.
If we denote by $X^{(1/2)}$ the normalization of $X$ inside
$k(X)^{1/2}$, then $X^{(1/2)}$ is abstractly isomorphic to $X$,
and the field extension $k(X)\subset k(X)^{1/2}$ 
induces a purely inseparable and finite morphism 
$F:X^{(1/2)}\to X$ of degree $4$, the $k$-linear Frobenius morphism.
Similarly, $k(\widetilde{X})\subset k(X)^{1/2}$ induces
a purely inseparable and finite morphism $\varpi:X^{(1/2)}\to\widetilde{X}$ 
of degree $2$ such that
$F=\pi\circ\varpi$.
Thus, we obtain a diagram
\begin{equation}
 \label{frobenius}
 \xymatrixcolsep{3pc}
 \xymatrix{
  X^{(1/2)} \ar[ddr]_F \ar[dr]|-{\varpi} \\
  & \widetilde{X} \ar[d]^\pi \ar[r]_\varphi & \PP^{1+C^2} \\
  & X  }
\end{equation}
The composition $\varphi\circ\varpi$ corresponds to a linear
subsystem of $|2C|$ (here, we identify $X$ with $X^{(1/2)}$).
Both, $\varphi$ and $\varpi$ are morphisms, and we have
$2\deg\varphi=\deg(\varphi\circ\varpi)$, as well
as $(2C)^2=4C^2$.
As before, we find $\deg\varphi\leq2$,
this time by arguing on $X^{(1/2)}$.

In order to show $\deg\varphi=1$ (now again, for arbitrary $\pi$ and 
$p$), we argue as in the proof of \cite[Lemma 4.4.3]{Cossec Models}.
Seeking a contradiction, we assume $\deg\varphi\neq1$.
Then, $\deg\varphi=2$ and the image $\varphi(\widetilde{X})$
is an integral surface of degree
$C^2$ in $\PP^{1+C^2}$, that is, a surface of minimal degree.
These surfaces have been explicitly classified by
del~Pezzo, but we refer to \cite{eh} for a characteristic-free discussion.

Now, the morphism $\pi$ is a torsor under a finite flat group
scheme $G$, which is of length $2$ over $k$.
Since the quotient of $\widetilde{X}$ by $G$ is isomorphic
to $X$ and not isomorphic to 
$\varphi(\widetilde{X})$ it follows that 
the $G$-action on $\widetilde{X}$ induces a non-trivial
$G$-action on $\PP(H^0(\widetilde{X},\pi^\ast\OO_X(C)))$ 
and $\varphi(\widetilde{X})$.
As already seen above, we can
write global sections of $\pi^\ast\OO_X(C)$ as
\begin{equation}
\label{seq:g-mod}
0\,\to\,H^0(X,\OO_X(C))\,\to\,
H^0(\widetilde{X},\pi^\ast\OO_X(C))\,\stackrel{{\rm pr}}{\longrightarrow}\,
H^0(X,\omega_X(C))\,\to\,0\,
\end{equation}
and consider it as sequence of $G$-modules.
It is not difficult to see (we will recall and generalize
this fact in Section \ref{subsec:explicit}) that 
$H^0(X,\OO_X(C))$ is the ${\rm id}$-eigenspace for
the $G$-action on $H^0(\widetilde{X},\pi^\ast\OO_X(C))$
and that $G$ acts via the determinant of its regular 
representation on $H^0(X,\omega_X(C))$.
%From the discussion at the beginning of this section
%it follows that $H^0(X,\pi_\ast\pi^\ast{\cal L})$, as
%a $G$-module, decomposes into the direct sum of three 
%regular representations of $G$, see also
%Proposition \ref{three quadrics 2a} below.

We set $\PP_+:=\PP(H^0(X,\OO_X(C)))$.
In case $G$ is linearly reductive, that is,
if $p\neq2$ or if $p=2$ and $G\iso\mu_2$, then the $G$-action has a second
eigenspace, providing us with a splitting of (\ref{seq:g-mod}), and which
we can identify with $H^0(X,\omega_X(C))$.
We denote by $\PP_-$ its projectivization and set
$\PP_-:=\emptyset$ in case $G$ is not linearly reductive.
Clearly, if a point in $\PP(H^0(\widetilde{X},\pi^\ast\OO_X(C)))$ 
is fixed under the $G$-action 
(in the scheme-theoretic sense) then it lies in $\PP_+$ or
$\PP_-$.
%, cf. also \cite[Lemma 4.4.3]{Cossec Models}.

For $v\in H^0(X,\omega_X(C))$,
the hyperplane $\PP_v:=\PP({\rm pr}^{-1}(v))$ is $G$-stable
and contains $\PP_+$.
For generic $v$, the intersection 
$\Delta:=\PP_v\cap\varphi(\widetilde{X})$ is an irreducible and
non-degenerate curve inside $\PP_v\iso\PP^{C^2}$.
Since $\Delta$ is of degree $C^2-1$ in a  
$C^2$-dimensional projective space,
it is a rational normal curve and
in particular, smooth and rational.
Since $\Delta$ is isomorphic to $\PP^1$ and 
equipped with a non-trivial $G$-action, its
fixed point scheme has length $2$.
% which is supported in two distinct points if $p\neq2$.

In particular, $\varphi(\widetilde{X})$ contains points that are
fixed under $G$ and so, its intersection with $\PP_+$
or $\PP_-$ is non-empty.
On the other hand,
$$
\PP_+\,\cap\,\varphi(\widetilde{X})\,=\,
\bigcap_{s\in \pi^\ast H^0(X,\OO_X(C))^\vee} \{s=0\}\,\cap\,\varphi(\widetilde{X}),
$$
and similarly for $\PP_-\cap\varphi(\widetilde{X})$ 
and $s\in\pi^\ast H^0(X,\omega_X(C))^\vee$.
This implies that $\OO_X(C)$ or $\omega_X(C)$ is not globally generated,
a contradiction.
Thus, $\deg\varphi=1$.

By \cite[Th\'eor\`eme I.6.10]{Jouanolou}, 
a generic
Cartier divisor in $|\pi^\ast\OO_X(C)|$ is irreducible.
The same theorem, applied to the open and dense subset of $\widetilde{X}$
where $\varphi$ is an isomorphism and where $\widetilde{X}$ is smooth,
shows that
a generic Cartier divisor is generically reduced.
Now, Cartier divisors on Gorenstein schemes are Gorenstein 
and in particular Cohen--Macaulay.
Thus, a generic Cartier divisor in $|\pi^\ast\OO_X(C)|$ is 
irreducible, generically reduced and Cohen--Macaulay.
Being generically reduced and Cohen--Macaulay, $C$ is reduced.
In particular, $C$ is integral.

Using the adjunction formula and $H^1(\widetilde{X},\OO_{\widetilde{X}})=0$,
we conclude that $\varphi$ induces on $D$ the 
morphism associated to the canonical sheaf $\omega_D$.
Thus, $\varphi$ being birational, the generic $D$
is not hyperelliptic.
(Hyperelliptic in the non-smooth case simply means that there exists
a finite morphism of degree $2$ onto $\PP^1$, see \cite{Schreyer}.)
% or \cite[\S1]{Liedtke Noether}.
\qed\medskip

Using Saint-Donat's analysis \cite{Saint-Donat}
of linear systems on K3 surfaces,
Cossec \cite{Cossec Models} has shown that $\varphi(\widetilde{X})$
in characteristic $\neq2$ is cut out by quadrics:
%Let us first slightly extend this result: 

\begin{Theorem}[Cossec$+\varepsilon$]
 \label{sdx}
 Under the assumptions of Theorem \ref{sd1},
 the image $\varphi(\widetilde{X})$ in $\PP^{1+C^2}$ is 
 projectively normal and cut out by quadrics, whenever
 \begin{enumerate}
  \item ${\rm char}(k)\neq2$, or
  \item ${\rm char}(k)=2$ and $X$ is a singular Enriques surface.
 \end{enumerate}
\end{Theorem}

\prf
In characteristic $\neq2$, this is shown in \cite[Section 8]{Cossec Models}.
%or \cite[Theorem 4.9.2]{Cossec; Dolgachev}.

If $X$ is a singular Enriques surface in characteristic $2$,
then $\widetilde{X}$ is smooth, $\varphi$ is birational
and $\varphi(\widetilde{X})$ has only isolated singularities.
But then, a generic divisor $D\in|\pi^\ast\OO_X(C)|$ is smooth
and the whole analysis in \cite[Section 7]{Saint-Donat}
remains valid also in characteristic $2$.
Thus, \cite[Theorem 7.2]{Saint-Donat} and
\cite[Lemma 8.1.2]{Cossec Models} show that $\varphi(\widetilde{X})$
is projectively normal and cut out by quadrics.
\qed\medskip

It is plausible that $\varphi(\widetilde{X})$ 
is always cut out by quadrics -- this would follow from
numerical $4$-connectedness of $\pi^\ast\OO_X(C)$ on
$\widetilde{X}$ (we refer to \cite[Section 3]{cfhr} for a discussion of
this notion for singular varieties).
However, we only establish this in 
the special case, in which we are interested in later
on:

\begin{Proposition}
 \label{special case}
 In addition to the assumptions of Theorem \ref{sd1}, assume 
 that $C^2=4$ and $\Phi(C)=2$. 
 Then, $\varphi(\widetilde{X})\subset\PP^5$ is 
 projectively normal and cut out by quadrics.
\end{Proposition}

\prf
In view of Theorem \ref{sdx} it suffices to treat the case
where $X$ is classical or supersingular in characteristic $2$.
Let us note that the proof we now present works in general,
but that we shall only give details in the case
where $\pi:\widetilde{X}\to X$ is purely inseparable,
that is, $X$ is classical or supersingular in $p=2$.

We have to show that the graded ring associated to
$\pi^\ast\OO_X(C)$ on $\widetilde{X}$ is generated in degree $1$
with relations in degree $2$ only, that is, a Koszul algebra.

By Theorem \ref{sd1}, a generic Cartier divisor 
$D\in|\pi^\ast\OO_X(C)|$ is an integral and non-hyperelliptic
Gorenstein curve of arithmetic genus $p_a(D)=5$.
For $n\geq1$, we consider the following
short exact sequences on $\widetilde{X}$:
$$
0\,\to\,\pi^\ast\OO_X((n-1)C)\,\to\,\pi^\ast\OO_X(nC)
\,\to\,\omega_D^{\otimes n}\,\to\,0\,.
$$
Pushing $\pi^\ast\OO_X((n-1)C)$ forward to $X$ and using 
% the vanishing results 
\cite[Theorem 1.5.1]{Cossec; Dolgachev},
we conclude $h^1(\widetilde{X},\pi^\ast\OO_X((n-1)C))=0$
for $n\geq1$.
Thus, as explained in the proof of part (ii) of
\cite[Theorem 6.1]{Saint-Donat}, to prove
our assertion, it suffices to show that the canonical 
ring of $D$ is a Koszul algebra.

Before proceeding, we study genus-one half-pencils
on $X$ more closely:
since $\Phi(C)=2$, there exists a genus-one pencil $E$
on $X$ with $C\cdot E=4$.
Moreover, let $E'$ be a genus-one curve
with $|2E'|=|E|$, that is, a half-pencil.
We now claim:
\begin{enumerate}
 \item $\pi^\ast\OO_X(E')$ is globally generated with $h^0=2$ and $h^1=0$, and
 \item $\pi^\ast\OO_X(C-E')$ also satisfies $h^0=2$ and $h^1=0$. It is 
   globally generated outside $\pi^{-1}(R)$, where
   $R$ is a (possibly empty) union of ${ADE}$-curves.
\end{enumerate}
We only deal with the case that $\pi$ is inseparable in characteristic $2$ 
and leave the remaining (and easier) cases to the reader:
the assertions about $h^0$ and $h^1$ of $\pi^\ast\OO_X(E')$
follow from pushing it 
down to $X$ and then using $h^0(X,\OO_X(E'))=h^0(X,\omega_X(E'))=1$, 
as well as $h^1(X,\OO_X(E'))=h^1(X,\omega_X(E'))=0$.
Now, we use diagram (\ref{frobenius}):
the linear system $\varpi^\ast|\pi^\ast\OO_X(E')|$
is a linear subsystem of $|\varpi^\ast\pi^\ast\OO_X(E')|=|F^\ast E'|=|2E'|=|E|$.
Since both satisfy $h^0=2$, they are equal and so $\pi^\ast\OO_X(E')$
is globally generated since $\OO_X(E)$ is.
Let us adjust the proof of \cite[Theorem 5.3.6]{Cossec Models} to our
situation:
by Riemann--Roch, there exists an effective divisor $D$ such that 
$E'+D\in|C|$, $E'D=2$, and $D^2=0$.
Moreover, there exists a divisor $E''$ of canonical type such that
$D=E''+R$ with $R\geq0$.
Since $\Phi(C)=2$, we have $CE''\geq2$ and if equality holds
then $E''$ is a genus-one half-pencil.
Thus, $4=C^2\geq CE'+CE''\geq4$, and we conclude
$CE''=2$ and $CR=0$.
In particular, if non-empty, $R$ is a union of ${ADE}$-curves.
The remaining assertions now follow as before, establishing our
two claims.
% problem with \cite[Proposition 3.6.2]{Cossec; Dolgachev} - I do not know
% whether C is 'hyperelliptic', same problem with 
% \cite[Theorem 5.3.6]{Cossec Models}.

Next, let us show that $\varphi(D)$ possesses a simple trisecant:
first, we choose a generic Cartier divisor $H\in|\pi^\ast\OO_X(E')|$.
Since $h^1(\widetilde{X},\pi^\ast\OO_X(C-E'))=0$, we conclude that
$H^0(\widetilde{X},\pi^\ast\OO_X(C))$ surjects onto
$H^0(G,\pi^\ast\OO_X(C)|_H)$.
Using $\deg\pi^\ast\OO_X(C)|_H=4$, we see that
$\varphi$ embeds $H$ as a quartic curve into some
$\PP^3$, which is easily seen to be the complete 
intersection of two quadrics.
% see the analysis in Eisenbud's ``The Geometry of Syzygies'', Chapter 6D
Thus, a generic hyperplane in $\PP^5$ intersects $H$
in $4$ points in uniform position.
This hyperplane cuts out on $\varphi(\widetilde{X})$ 
an integral curve $\varphi(D)$, where $D\in|\pi^\ast\OO_X(C)|$,
having the stated simple trisecant. 

Having established this trisecant, and noting that
$p_a(D)-2=3$, \cite[Theorem 1.2]{Schreyer} and 
\cite[Corollary 1.3]{Schreyer} show that the 
canonical ring of $D$ is generated in degree $1$ 
and has relations in degree $\leq3$.
Our proposition is proved once we show that no relations
in degree $3$ are needed.

Suppose that relations in degree $3$ are needed.
Then, $\varphi(D)$ is contained in an irreducible 
surface $S$ of degree $p_a(D)-2=3$ by \cite[Theorem 3.1]{Schreyer}.
By the classification of surfaces of minimal degree
\cite{eh} together with $p_a(D)=5$ we find that $S$ is ruled.
Moreover, a generic ruling of $S$ intersects $\varphi(D)$
in three distinct smooth points.
Thus,
$D$ possesses a globally generated invertible sheaf $\cal M$ 
of degree $3$ with $h^0=2$ (a $g_3^1$ in classical terminology) 
and $D$ is trigonal. 
Now, we consider ${\cal L}:=\pi^\ast\OO_X(E')|_D$. 
Then $\pi^\ast\OO_X(C-E')|_D\iso\omega_D\otimes{\cal L}^\vee$
and ${\cal L}$ and $\omega_D\otimes{\cal L}^\vee$
are invertible sheaves of degree $4=\frac{1}{2}\deg\omega_D$.
Taking cohomology in 
$$
0\,\to\,\pi^\ast\OO_X(E'-C)\,\to\,\pi^\ast\OO_X(E')\,\to\,{\cal L}\,\to\,0,
$$
we obtain $h^0(D,{\cal L})\geq2$.
Moreover, $H^0(\widetilde{X},\pi^\ast\OO_X(E'))$
and $H^0(\widetilde{X},\pi^\ast\OO_X(C-E'))$
inject into $H^0(D,{\cal L})$ and 
$H^0(D,\omega_D\otimes{\cal L}^\vee)$, respectively.
In particular, ${\cal L}$ is globally generated since 
$\pi^\ast\OO_X(E')$ is.
Choosing $D$ to be generic, we may assume that $D$ does not
intersect $\pi^{-1}(R)$ and thus,
$\omega_D\otimes{\cal L}^\vee$ is globally generated.
Since $h^1(D,{\cal L})\neq0$, Clifford's inequality implies
$h^0(D,{\cal L})\leq3$ and equality could only happen
if $D$ was hyperelliptic.
Thus, $h^0(D,{\cal L})=h^0(D,\omega_D\otimes{\cal L}^\vee)=2$.
This is enough to show that $D$ is not trigonal:
by \cite[Excercise III.B-5]{acgh}, which works for
integral Gorenstein curves in arbitrary characteristic,
the invertible sheaf $\cal M$ making $D$ trigonal %(the ``$g_3^1$'')
would have to be a subsheaf of ${\cal L}$ or 
$\omega_D\otimes{\cal L}^\vee$, which is absurd.
\qed\medskip

Although we will not need this result in the sequel, let us
complete the picture by discussing linear systems on $\widetilde{X}$
arising from curves $C$ on $X$ with $\Phi=1$:

\begin{Theorem}
 \label{sd3}
% Let $X$ be an Enriques surface and
 Let $C$ be an irreducible curve with $C^2>0$ and $\Phi(C)=1$
 on an Enriques surface $X$.
 Then
 \begin{enumerate}
  \item the invertible sheaf $\pi^\ast\OO_X(C)$ on $\widetilde{X}$
    is globally generated,
  \item $|\pi^\ast\OO_X(C)|$ gives rise to a morphism 
    $$
   \varphi\,:\,\widetilde{X}\,\to\,\PP^{(1+C^2)},
   $$
   which is generically of degree $2$ onto
   a surface of minimal degree $C^2$,
  \item the image $\varphi(\widetilde{X})$ is cut out by quadrics.
\end{enumerate}
\end{Theorem}

\prf
As in the proof of Theorem \ref{sd1} we find
$h^0(\widetilde{X},\pi^\ast\OO_X(C))=2+C^2$.
Again, $|C|$ has no fixed component
%(by the same argument as in the proof of Theorem \ref{sd1}),
and so $\pi^\ast\OO_X(C)$
is globally generated outside a finite (possibly empty) 
set of points.

Let us first assume $C^2\geq4$.
In this case, $\varphi(\widetilde{X})$ is a surface since the 
image of the rational map associated to $|C|$ is a surface
\cite[Theorem 4.5.1]{Cossec; Dolgachev}.

Seeking a contradiction, we assume that $\varphi$ is birational.
As in the proof of Theorem \ref{sd1}, we conclude
that a generic Cartier divisor $D\in|\pi^\ast\OO_X(C)|$ 
is an integral Gorenstein curve.
Since $\Phi(C)=1$, there exists
a genus-one half-pencil $E'$ on $X$ such that $C\cdot E'=1$.
Then ${\cal L}:=\pi^\ast\OO_X(E')|_D$ satisfies
$\deg{\cal L}=2$ and taking cohomology in
$$
0\,\to\,\pi^\ast\OO_X(E'-C)\,\to\,\pi^\ast\OO_X(E')\,\to\,{\cal L}\,\to\,0
$$
we find $h^0(D,{\cal L})\geq2$.
Since $p_a(D)\geq5$, Riemann--Roch implies $h^1(D,{\cal L})\neq0$.
But then, Clifford's inequality $h^0(D,{\cal L})\leq2$ is in fact 
an equality, which implies that $D$ is hyperelliptic.
In the proof of Theorem \ref{sd1} we have seen that $\varphi$
restricted to $D$ induces $|\omega_D|$, 
which contradicts the fact that $\varphi$ is birational.
Thus, $\deg\varphi\geq2$ and since
$\varphi(\widetilde{X})$ is a non-degenerate integral surface
in $\PP^{1+C^2}$, we conclude
$$
2C^2\,\leq\,\deg\varphi\cdot C^2\,\leq\,
\deg\varphi\cdot\deg\varphi(\widetilde{X}).
$$

On the other hand, $\pi^\ast\OO_X(C)$ is globally generated outside
a finite set of points and so we find
$$
\deg\varphi\cdot\deg\widetilde{X}\,\leq\,2C^2
$$
with equality if and only if $\pi^\ast\OO_X(C)$ is globally generated:
this is clear if $\pi$ is \'etale, because then $\widetilde{X}$ is smooth.
If $\pi$ is inseparable, we consider $\varphi\circ\varpi$ in the
diagram (\ref{frobenius}) and obtain the same result by arguing on
$X^{(1/2)}$.

Putting these inequalities together, we find that $\pi^\ast\OO_X(C)$
is globally generated, $\deg\varphi=2$ and 
$\deg\varphi(\widetilde{X})=C^2$.
In particular, $\varphi(\widetilde{X})$ is a surface of minimal 
degree and thus, cut out by quadrics \cite{eh}.

It remains to deal with the case $C^2=2$.
Then, $\varphi$ is a possibly rational map to $\PP^3$.
By contradiction, assume that $\varphi(\widetilde{X})$ is a curve.
A generic $H\in|\pi^\ast\OO_X(E')|$, where $E'$ is a genus-one
half-pencil with $C\cdot E'=1$, is an integral curve with $p_a=1$.
We find $\deg\pi^\ast\OO_X(C)|_H=2$, which implies
$h^0(H,\pi^\ast\OO_X(C)|_H)=2$ by Riemann--Roch and Clifford's
inequality.
This implies that $\varphi(H)$ is a linearly embedded $\PP^1\subset\PP^3$.
But then $\varphi(\widetilde{X})$ is equal to this $\PP^1$,
contradicting that $\varphi(\widetilde{X})$ linearly spans $\PP^3$.
Thus, $\varphi(\widetilde{X})$ is a surface and we conclude
as before.
\qed

\section{Complete intersections of three quadrics}
\label{sec:three quadrics}

For an Enriques surface $X$, we study in this section 
a certain birational and projective model of its canonical double cover
$\widetilde{X}$, which will turn out to be
a complete intersection of three quadrics in $\PP^5$.
This extends results of Cossec \cite[Section 8]{Cossec Models}
and Verra \cite[Theorem 5.1]{Verra} - in some form already
known to Enriques - to characteristic $2$, and rests on
Proposition  \ref{special case} from the previous section.
We explicitly describe the equations and
the action of the finite flat group scheme $G=\Pic^\tau(X)^D$
on $\widetilde{X}$.
In particular, we prove that {\em all}
Enriques surfaces in {\em any} characteristic
arise via the Bombieri--Mumford--Reid construction 
presented in \cite[\S3]{bm3}.

\subsection{Birational models of the canonical double cover}

In order to find projective models of $\widetilde{X}$,
we study linear systems $|\pi^\ast\OO_X(C)|$, 
where $C\subset X$ is an irreducible curve with $C^2>0$. 
By Theorem \ref{sd1} and Theorem \ref{sd3},
the associated morphism 
$\varphi:\widetilde{X}\to\PP^N$
is birational onto its image
if and only if $\Phi(C)\geq2$ and $C^2\geq4$.
In this case, the codimension of the image %$\varphi(\widetilde{X})$ 
is equal to $C^2-1$.
In this setup, models of smallest codimension
are of codimension $3$ in $\PP^5$.
Moreover,
by \cite[Lemma 3.6.1]{Cossec; Dolgachev}, irreducible
curves with $C^2<10$ satisfy $\Phi(C)\leq2$.
We are thus led to studying irreducible curves
with $C^2=4$ and $\Phi(C)=2$.

\begin{Theorem}
 \label{three quadrics}
 For every Enriques surface $X$ there exists a morphism 
 $\varphi:\widetilde{X}\to\PP^5$ that is birational onto
 its image.
 More precisely, there is a Cartesian diagram
% and an embedding 
 $$
  \xymatrixcolsep{2pc}
  \xymatrix{
     \varphi\ar@{}[r]|{:} & \widetilde{X} \ar@{->}[r]\ar@{->}[d]_\pi&%\ar@{}[dr]|{\Box} &
     \varphi(\widetilde{X})\ar@{->}[d]_{\pi'}\ar@{^{(}->}[r] &
     \PP^5\\
    & X \ar@{->}[r]^{\nu} & X'
 }
 $$
% $$\begin{array}{ccccccc}
%   \varphi &:&\widetilde{X}&\to&\varphi(\widetilde{X})&\subset&\PP^5\\
%   &&\downarrow_\pi&&\downarrow_{\pi'}\\
%   &&X&\stackrel{\nu}{\to}&X'
%  \end{array}$$
 such that
 \begin{enumerate}
%  \item $\varphi:\widetilde{X}\to\varphi(\widetilde{X})$ is a birational morphism,
  \item $\varphi(\widetilde{X})$ is a complete intersection of three quadrics,
  \item $\nu$ is a birational morphism and $X'$ has at worst Du~Val singularities,
  \item $\pi$ is a torsor under the finite flat group scheme $G:=(\Pic^\tau(X))^D$ that arises
    as pull-back from a $G$-torsor $\varphi(\widetilde{X})\to X'$, 
  \item the $G$-action on $\varphi(\widetilde{X})$ is induced by a linear $G$-action of the
    ambient $\PP^5$.
 \end{enumerate}
\end{Theorem}

\prf
By \cite[Chapter IV \S9]{Cossec; Dolgachev}, there
exists an irreducible curve $C\subset X$ with
$C^2=4$ and $\Phi(C)=2$.
% , but see also Proposition \ref{prop:polarization}.
We set ${\cal L}:=\OO_X(C)$ and then,
by Theorem \ref{sd1} and Proposition \ref{special case}, 
the invertible sheaf $\pi^*{\cal L}$ 
on $\widetilde{X}$ gives rise to a birational morphism,  
whose image $\varphi(\widetilde{X})$ is cut out by quadrics.
Taking cohomology in the short exact sequences
$$
0\,\to\,{\cal L}^{\otimes n}\,\to\,\pi_\ast\pi^\ast({\cal L}^{\otimes n})\,\to\,
\omega_X\otimes{\cal L}^{\otimes n}\,\to\,0,
$$
using Riemann--Roch on $X$ and the vanishing result 
\cite[Theorem 1.5.1]{Cossec; Dolgachev}, we find
 $h^0(\widetilde{X},\pi^*{\cal L})=6$ and 
$h^0(\widetilde{X},\pi^*{\cal L}^{\otimes 2})=18$.
%using the projection formula, (\ref{ses}) and the results of 
%\cite[Chapter I \S5]{Cossec; Dolgachev}.
Thus, there are three quadric relations and hence,
$\varphi(\widetilde{X})$ is a complete intersection of three
quadrics in $\PP^5$.

Next, the $G$-action on $\widetilde{X}$ induces a $G$-action
on $H^0(\widetilde{X},\pi^*{\cal L})$.
This gives rise to a linear $G$-action on $\PP^5$ extending the
$G$-action on $\varphi(\widetilde{X})$.

%In the proof of Proposition \ref{special case} we have seen that
%we may decompose $C$ as $E'+E''+R$, where $E'$, $E''$ are genus one 
%half-pencils (possibly equal), 
%and $R$ is a (possibly empty) union of $(-2)$-curves.
By \cite[Proposition 4.1.1]{Cossec; Dolgachev}, every irreducible curve that has 
zero-intersection with $C$ is a $(-2)$-curve.
Since $\OO_X(C)$ is globally generated, big and nef, 
% we conclude that 
$$
  \nu\,:\,X \,\longrightarrow\, X'\,:=\,\Proj \bigoplus_{n\geq0} H^0(X,\,{\cal L}^{\otimes n})
$$ 
is a birational morphism that contracts those  $(-2)$-curves 
having zero-intersection with $C$ and does not contract anything else.
In particular, $X'$ has at worst Du~Val singularities.
Thus, $H^1(X,\OO_X)\iso H^1(X',\OO_{X'})$ and $\omega_X$ is $2$-torsion
if and only if $\omega_{X'}$ is.
In particular, the canonical $G$-torsor $\pi$ arises as pull-back
from a $G$-torsor  $\widetilde{X}'\to X'$.

Since $X'$ has only Du~Val singularities, 
${\cal L}^{\otimes n}$ (resp. $\omega_X\otimes{\cal L}^{\otimes n}$) 
for $n\geq0$ and $\nu_\ast({\cal L}^{\otimes n})$
(resp. $\nu_\ast(\omega_X\otimes{\cal L}^{\otimes n})$) have isomorphic global sections.
Thus, the graded ring $\widetilde{R}_{\pi^*\cal L}$ of 
$(\widetilde{X},\pi^*{\cal L})$ is
isomorphic to the graded ring $\widetilde{R}'_{\pi'^*\nu_*{\cal L}}$ 
of $(\widetilde{X}',\pi'^*\nu_\ast{\cal L})$.
By Proposition \ref{special case}, we can identify
$(\Proj \widetilde{R}_{\pi^*\cal L},\OO(1))$ with
$\varphi(\widetilde{X})\subset\PP^5$.
On the other hand, $\nu_\ast{\cal L}$ is ample on $X'$
by the Nakai--Moisehzon criterion, and so,
$\nu_\ast\pi^*{\cal L}\iso\pi'^\ast\nu_\ast{\cal L}$
is ample on $\widetilde{X}'$.
Thus, $\widetilde{X}'$ is isomorphic to $\Proj \widetilde{R}'_{\pi^*\cal L}$.
\qed

\subsection{Cossec--Verra Polarizations}
Having just established a birational projective model
of $\widetilde{X}$, it is natural to ask for uniqueness
of $\varphi:\widetilde{X}\to\PP^5$, as well as to understand how 
close $\varphi$  is to being an isomorphism onto its image.
Since our previous result extends work of Cossec \cite[Section 8]{Cossec Models} 
and Verra \cite[Theorem 5.1]{Verra} 
to characteristic $2$, we define the following notion:

\begin{Definition}
 \label{def:polarization}
 A {\em Cossec--Verra polarization} on an Enriques surface $X$
 is an invertible sheaf ${\cal L}\in{\cal P}(X)$, where 
 $$
    {\cal P}(X)\,:=\,
   \left\{ 
      {\cal L}\in\Pic(X)\, \left|\,
            \begin{array}{l}
              \mbox{there exists an irreducible curve }C\\
              \mbox{with }C\in|{\cal L}|,\,C^2=4,\, \Phi(C)=2 
           \end{array}
   \right\} \right. .
 $$ 
 We denote the morphism $\widetilde{X}\to\PP^5$
 corresponding to $|\pi^\ast{\cal L}|$ by $\varphi_{\cal L}$.
% Moreover, if $X'$ has at worst Du~Val singularities and
% its minimal desingularization $\nu:X\to X'$ is an Enriques surface,
% an invertible sheaf ${\cal L}$ on $X'$ is called a Cossec--Verra
% polarization if $\nu^\ast{\cal L}$ on $X$ is. 
\end{Definition}

%Then, $\varphi$ is equal to $\varphi_{\cal L}$, where $\cal L$ is a
%Cossec--Verra polarization.
By Theorem \ref{three quadrics}, 
Cossec--Verra polarizations are big and nef, but need not be ample.
Thus, it might be more appropriate to talk about
Cossec--Verra quasi-polarizations, but for the sake of readability we have
decided not to do so.
To decide whether a Cossec--Verra polarization is ample, let us recall
that an irreducible curve on $X$ is a $(-2)$-curve, 
that is, has self-intersection $-2$,
if and only if it is smooth and rational.
Such curves are called {\em nodal} and we denote by
the set of all nodal curves by ${\cal R}(X)$.

\begin{Proposition}
 \label{prop:elementary}
 For 
 ${\cal L}\in{\cal P}(X)$ the following properties are equivalent:
 \begin{enumerate}
  \item $\deg{\cal L}|_\alpha>0$ for every $\alpha\in{\cal R}(X)$,
  \item ${\cal L}$ is ample,
  \item $\pi^\ast{\cal L}$ is very ample, and
  \item $\varphi_{\cal L}:\widetilde{X}\longrightarrow\varphi_{\cal L}(\widetilde{X})$ is
   an isomorphism.
 \end{enumerate}
 In general, the reduced exceptional locus of $\varphi_{\cal L}$ is the union of the reduced
 inverse images under $\pi$ of all those nodal curves $\alpha\in{\cal R}(X)$
 with $\deg{\cal L}|_\alpha=0$. 
\end{Proposition}

\prf
By \cite[Corollary 3.2.2]{Cossec; Dolgachev}, every effective divisor
is linearly equivalent to one that is the positive sum of
curves of arithmetic genus $1$ and $0$.
Since $\Phi(C)=2$, the intersection of ${\cal L}$ with curves of arithmetic genus $1$
is positive.
Thus, by the Nakai--Moishezon criterion for
ampleness, ${\cal L}$ is ample if and only if it has positive intersection with
every nodal curve.
This establishes $1\Leftrightarrow 2$.
From the proof of Theorem \ref{three quadrics} we get $1\Leftrightarrow 4$.
The equivalence $3\Leftrightarrow4$ is a tautology.

Finally, the last assertion follows from the proof of Theorem \ref{three quadrics}.
\qed\medskip

An Enriques surface $X$ is called 
{\em unnodal} if ${\cal R}(X)$ is empty.
Over the complex numbers, a generic Enriques surface is
unnodal \cite[Proposition 2.8]{Barth; Peters}.
We recall that a {\em Reye congruence} is the subset of 
the Gra\ss mannian 
$\GG(1,3)$ parametrizing lines in $\PP^3$ that lie on
at least two quadrics of a given generic $3$-dimensional family
of quadrics.
By \cite[Theorem 1]{Cossec Reye},
a generic nodal Enriques surfaces is a Reye congruence.
A rather extremal class of Enriques surfaces is called {\em extra-special},
which exists in characteristic $2$ only, and which is defined by the existence
of certain configurations of nodal curves (see \cite[Chapter III.5]{Cossec; Dolgachev}
for the precise definition).
This said, we establish a couple of facts about existence, uniqueness
and ampleness of Cossec--Verra polarizations:

\begin{Proposition}
 \label{prop:polarization}
 Every Enriques surface $X$ possesses at least one Cossec--Verra polarization,
 that is, ${\cal P}(X)\neq\emptyset$.
 Moreover,
 \begin{enumerate}
 \item if $X$ is unnodal, then ${\cal P}(X)$ is infinite and 
   every ${\cal L}\in{\cal P}(X)$ is ample, 
 \item if $X$ is a generic Reye congruence and $p\neq2$, then
   there exists an ample Cossec--Verra polarization and 
   ${\cal P}(X)$ is infinite, 
 \item if $X$ is an extra special Enriques surface in characteristic $p=2$,
    then ${\cal P}(X)$ is finite and no 
    Cossec--Verra polarization is ample.
    Moreover, there exist nodal curves on $X$, whose inverse images on $\widetilde{X}$
    are contracted by $\varphi_{\cal L}$ for every ${\cal L}\in{\cal P}(X)$.
% \item if $X$ is $\widetilde{{\sf E}}_8$-extra special then there exists 
%    only one Cossec--Verra polarization ${\cal L}$.
%    Moreover, $\widetilde{X}$ is non-normal and the corresponding
%    $\varphi_{\cal L}$ partially
%    contracts the non-normal locus.
 \end{enumerate}
 Over the complex numbers, the set ${\cal P}(X)$ modulo ${\rm Aut}(X)$ is finite.
 More precisely, we have
 $$ 
       | {\cal P}(X) /  {\rm Aut}(X) | \,=\, 252,960 
 $$
 if $X$ is a generic complex Enriques surface.
\end{Proposition}

\begin{Remark}
 \label{rem:extra special not ample}
 Extra special Enriques surfaces exist in characteristic $2$ only and
 are discussed in \cite[Chapter III.5]{Cossec; Dolgachev}.
 For example, an 
 $\widetilde{{E}}_8$-extra special Enriques surface
 possesses only one Cossec--Verra polarization,
 $\widetilde{X}$ is non-normal and the unique
 $\varphi_{\cal L}$ partially contracts the non-normal locus.
 We shall see in Section \ref{sec:final} below that this is 
 closely related to obstructions of the deformation functor, as well
 as to the exceptional Enriques surfaces studied in \cite{Ekedahl}.  
\end{Remark}

\prf
Let us recall, for example, from \cite[Chapter II.5]{Cossec; Dolgachev}
that the Enriques lattice $\EE$, that is,
the N\'eron--Severi group of an Enriques surface
modulo torsion, is isometric to the hyperbolic lattice
corresponding to the Dynkin diagram ${T}_{2,3,7}$:
$$
 \xymatrixcolsep{2pc}
 \xymatrix{
    \ar@{}[d]|{\alpha_1} &
    \ar@{}[d]|{\alpha_2} &
    \ar@{}[d]|{\alpha_3} &
    \ar@{}[d]|{\alpha_4} &
    \ar@{}[d]|{\alpha_5} &
    \ar@{}[d]|{\alpha_6} &
    \ar@{}[d]|{\alpha_7} &
    \ar@{}[d]|{\alpha_8} &
    \ar@{}[d]|{\alpha_9} \\
    \circ\ar@{-}[r] &  
    \circ\ar@{-}[r] &  
    \circ\ar@{-}[r]\ar@{-}[d] & 
    \circ\ar@{-}[r] &
    \circ\ar@{-}[r] &  
    \circ\ar@{-}[r] &  
    \circ\ar@{-}[r] &  
    \circ\ar@{-}[r] &  
    \circ \\
    & & \circ \\
    & & \ar@{}[u]|{\alpha_0} \\
  }
$$
Let $W_\EE$ be the Weyl group with respect to all roots of $\EE$ 
and let $W_X$ be the Weyl group with respect to the nodal curves ${\cal R}(X)$.
We define the {\em nodal chamber} of $X$ to be
$C_X:=\{x\in\EE\,|\, x\cdot\alpha\geq0,\forall\alpha\in{\cal R}(X)\}$.
As explained in \cite[Chapter III.2]{Cossec; Dolgachev},
$C_X$ is a fundamental domain of $V_X:=\{ x\in\EE\,|\,x^2\geq0\}$
for the $W_X$-action.
Moreover, let $V_X^+$ be the connected component of 
$V_X$ containing the class of an ample divisor and 
set $C_X^+:=V_X^+\cap C_X$.
%Next, $\EE$ is unimodular, and so we may and will identify it with its
%dual lattice $\EE^\vee$.

Let $\omega_1:=\alpha_1^\vee\in\EE^\vee$ be the fundamental weight
of the root $\alpha_1\in\EE$.
It follows from \cite[Corollary 2.5.7]{Cossec; Dolgachev} that
elements of ${\cal P}(X)$ correspond to those elements of the orbit
${\rm Isom}(\EE)\cdot\omega_1$ that lie in $C_X^+$.
In particular, ${\cal P}(X)$ is not empty.
The stabilizer ${\rm Stab}(\omega_1)$ is the Weyl group
corresponding to the Dynkin diagram ${T}_{2,3,7}$ with
vertex $\alpha_1$ removed, which is of type ${D}_9$.
In particular, ${\rm Stab}(\omega_1)$ is finite.

If $X$ is unnodal then $W_X$ is trivial, $C_X^+=V_X^+$ and so
${\cal P}(X)$ corresponds to the cosets of 
$W_\EE={\rm Isom}(\EE)/\{\pm{\rm id}\}$
modulo ${\rm Stab}(\omega_1)$, which is infinite.
Since $X$ is unnodal, every 
${\cal L}\in{\cal P}(X)$ is ample by Proposition \ref{prop:elementary}.
Moreover, from ${\rm Stab}(\omega_1)\iso W({D}_9)$,
we infer $W({D}_9)/W({D}_9)(2)\iso(\ZZ/2\ZZ)^8\rtimes\Sym_9$, see 
\cite[Proposition 2.8.4]{Cossec; Dolgachev}.
By \cite[Theorem (3.4)]{Barth; Peters}, a generic complex Enriques surface satisfies
${\rm Aut}(X)\iso W_\EE(2)$ and thus 
$W_\EE/W_\EE(2)\iso  O^+(10,\FF_2)$ by \cite[Theorem 2.9.1]{Cossec; Dolgachev}.
This identifies ${\cal P}(X)$ modulo ${\rm Aut}(X)$ with 
$W_\EE/W_\EE(2)$ modulo $W({D}_9)/W({D}_9)(2)$, which is of order
$$
   \frac{|O^+(10,\FF_2)|}{|(\ZZ/2\ZZ)^8\rtimes\Sym_9|}\,=\,
   \frac{2^{20}\cdot3^5\cdot 5^2\cdot7\cdot17\cdot31}{2^{15}\cdot 3^4\cdot 5\cdot 7}\,=\,
   2^5\cdot3\cdot5\cdot17\cdot31\,=\, 252,960   \,.
$$

Next, let $X$ be a generic Reye congruence.
It contains $10$ genus-one half-pencils $F_i$ 
and $10$ nodal curves $D_i$ such that $F_i F_j=1$ and
$D_i D_j=2$ for $i\neq j$, see \cite[Lemma 3.2.1]{Cossec Reye}.
It follows from the proof of \cite[Proposition 3.2.5]{Cossec Reye},
that the invertible sheaf corresponding to
$C_{ij}:=F_i+\frac{1}{2}(D_i+D_j)$ for $i\neq j$ belongs to ${\cal P}(X)$.
If $X$ is a generic Reye congruence then the genus-one fibrations
$|2F_i|$ have no reducible fibers by 
the remark after \cite[Proposition 3.2.4]{Cossec Reye}.
From this it is not difficult to compute that every
nodal curve intersects $C_{ij}$ positively, that is,
the corresponding $\varphi_{\cal L}$ is an isomorphism.
It follows from \cite[Theorem 1]{Cossec Dolgachev Automorphisms}, that
automorphism groups of generic nodal Enriques surfaces
in characteristic $p\neq2$ are infinite.
Since this group acts on $C_X^+$ and since ${\rm Stab}(\omega_1)$
is finite, we conclude that ${\cal P}(X)$ is infinite.

If $X$ is $\widetilde{{E}}_8$-extra special,
then $W_X={\rm Isom}(\EE)/\langle\pm{\rm id}\rangle$, which implies
that it contains only one genus-one fibration $|2E|$ and that
${\cal P}(X)$ consists of only one element $\cal L$.
Then we use \cite[Proposition 3.6.2]{Cossec; Dolgachev} to
see that $|{\cal L}|=|2E+2R_1+...+2R_7+R_8+R_{10}|$
(notation as in case $4$ of \cite[page 185]{Cossec; Dolgachev})
from which we read off that $R_1,...,R_8$ and $R_{10}$
have zero-intersection with ${\cal L}$.
In the other extra special cases, the genus-one fibrations
are described in \cite[Chapter III.5]{Cossec; Dolgachev}
and applying \cite[Proposition 3.6.2]{Cossec; Dolgachev}
to a divisor class $|C|$ with $C^2=4$, $\Phi(C)=2$ we end 
up with a finite list of possibilities of how to write
$|C|$ in terms of genus-one fibrations.
First, this shows that ${\cal P}(X)$ is finite.
Second, in these explicit lists we can always find nodal
curves that have zero-intersection with $C$ for
any choice of $C$ and any decomposition into genus-one
pencils.
We leave the lengthy, yet straight forward 
details to the reader.

Finally, over the complex numbers, 
the subgroup of ${\rm Isom}(\EE)$ generated by $W_X$
and ${\rm Aut}(X)$ is of finite index 
\cite{Dolgachev Automorphisms}.
In particular, there are only finitely many orbits
of ${\rm Isom}(\EE)\cdot \omega_1$ modulo $W_X$
(needed to move the vector into $C_X^+$) and modulo
${\rm Aut}(X)$. 
Thus, ${\cal P}(X)$ modulo ${\rm Aut}(X)$
is finite.
\qed

\subsection{Explicit description and equations}
\label{subsec:explicit}
We end this section by giving an explicit description of the canonical double
cover, its $G$-action, and the 
$G$-invariant quadrics cutting out $\varphi_{\cal L}(\widetilde{X})$
in $\PP^5$.

For later use, let us slightly enlarge our setup: we let $R$ be a complete,
local, and Noetherian ring with residue field $k$.
We let ${\cal X}\to\Spf R$ be a formal family of Enriques surfaces
with special fiber $X\to \Spec k$.
We also assume that there is an invertible sheaf $\overline{{\cal L}}$
on $\cal X$ that restricts to a Cossec--Verra polarization on the special fiber.
Since we don't know yet that this is always true 
(see Proposition \ref{prop:picard scheme flat} below),
we also assume $\Pic^\tau({\cal X}/R)$ exists and is a finite flat
group scheme of length $2$ over $\Spf R$.

We start by describing the canonical double
cover:  since we assumed $\Pic^\tau({\cal X}/R)$ to be
a finite and flat group scheme of length $2$, it gives rise to a torsor
$$
  \pi\,:\,\widetilde{{\cal X}}\,\longrightarrow\,{\cal X}
$$
under its Cartier dual $(\Pic^\tau({\cal X}/R))^D$, see
\cite[Proposition (6.2.1)]{Raynaud}.

More precisely, let ${\cal P}$ be the Poincar\'e invertible
sheaf on ${\cal X}\times\Pic^\tau({\cal X}/R)$.
By its universal property, there exists a 
morphism $\psi:\Pic^\tau({\cal X}/R)\to\Pic^\tau({\cal X}/R)$
such that ${\cal P}\otimes{\cal P}\iso({\rm id}\times\psi)^*{\cal P}$.
Clearly, $\psi=\mu\circ\Delta$, where $\Delta$ is the diagonal
and $\mu$ is the multiplication map of 
$\Pic^\tau({\cal X}/R)$.
Dualizing, we obtain an $\OO_{\cal X}$-algebra structure
on ${\cal P}^\vee$.
Dualizing the %$\OO_{\cal X}\otimes\OO_{\cal G}$-module
multiplication map
${\cal P}\otimes(\OO_{\cal X}\otimes\OO_{{\Pic^\tau({\cal X}/R)}})\to{\cal P}$,
we obtain a $\OO_{\cal X}\otimes\OO_{{\Pic^\tau({\cal X}/R)}^D}$-comodule
structure on ${\cal P}^\vee$.
Putting these observations together, we obtain the following
description of the canonical double cover
$$
   \pi\,:\,
   \widetilde{{\cal X}}\,\iso\,\BigSpec {\cal P}^\vee \,\longrightarrow\,{\cal X}
$$
together with its $(\Pic^\tau({\cal X}/R))^D$-action.
In particular, this group scheme acts via its regular representation.

By the Tate--Oort classification \cite[Theorem 2]{Tate; Oort},
a finite flat group scheme of length $2$ over $R$ is isomorphic to
$$
 {\cal G}_{a,b} \,:=\,\Spec R[x]/(x^2-ax) \,\mbox{ \quad with \quad } ab=2
$$
and the comultiplication is given by
$x\mapsto x\otimes1+1\otimes x-bx\otimes x$.
Note that $a,b$ are not unique since we have 
${\cal G}_{a,b}\iso{\cal G}_{ar,br^{-1}}$ for every unit $r\in R^\times$.
Also, the Cartier dual group scheme of ${\cal G}_{a,b}$ is isomorphic to
${\cal G}_{b,a}$.

Having recalled these facts from \cite{Tate; Oort}, 
there exist $a,b\in R$ with $ab=2$ such that 
$\Pic^\tau({\cal X}/R)\iso {\cal G}_{b,a}$, and then,
$\Pic^\tau({\cal X}/R)^D\iso{\cal G}_{a,b}$.
We obtain the {\em regular representation} of ${\cal G}_{a,b}$ by
assigning to every $R$-algebra $S$ the group homomorphism 
$$\begin{array}{ccccc}
    \rho_{\rm reg} &:& {\cal G}_{a,b}(S)\,=\,\{s\in S\,|\,s^2=as \}&\to&{\rm GL}_2(S) \\
    && s &\mapsto& \left(\begin{array}{cc}
             1 & s \\
             0 & 1-bs
            \end{array}\right)
   \end{array}
$$
see also \cite[\S3]{bm3}.

\begin{Lemma}[Bombieri--Mumford]
 \label{group action}
 We set $T:=R[x_1,y_1,...,x_n,y_n]$ and assume that 
 ${\cal G}_{a,b}$ acts on each pair $x_i,y_i$  
 via $\rho_{\rm reg}$. 
 Then, the following quadrics are ${\cal G}_{a,b}$-invariant
 $$
      x_ix_j,\mbox{ \qquad } y^2_i-a\,x_iy_i,\mbox{ \qquad }
      x_iy_j+y_ix_j+b\, y_iy_j\,.
 $$
 Moreover, the ${\cal G}_{a,b}$-invariants of $T$ in even degree
 are generated by these invariant quadrics.
\end{Lemma}

\prf
\cite[p.222]{bm3}.
\qed\medskip

% reason: these invariants have to be invariant for EVERY R-algebra S,
% i.e., nilpotent and other zero divisors in R do not give rise
% to 'special' or 'extra' invariants - thus, have to solve it UNIVERSALLY 

With these preparations, we now show that ${\cal G}_{a,b}$-invariant
quadrics as in the previous lemma cut out the image of ${\cal X}$
inside $\PP^5$:

\begin{Proposition}
  \label{three quadrics explicit}
  Let $f:{\cal X}\to\Spf R$ be a formal family of Enriques surfaces,
  together with $\overline{\cal L}\in\Pic({\cal X}/R)$ that restricts
  to a Cossec--Verra polarization on the special fiber.
  Also, assume that $\Pic^\tau({\cal X}/R)\iso{\cal G}_{b,a}$ 
  for some $a,b\in R$ with $ab=2$.
  Then, there exists a linear ${\cal G}_{a,b}$-action on $\PP^5_R$,
  such that $\pi^\ast\overline{\cal L}$ defines a
  ${\cal G}_{a,b}$-equivariant morphism
  $$
  \varphi_{{\overline{\cal L}}} \,:\,\widetilde{{\cal X}}\,\longrightarrow\,\PP_R^5
  $$
  over $\Spf R$, whose image is the complete intersection of three quadrics.
  More precisely, there exist coordinates
  $x_1,x_2,x_3,y_1,y_2,y_3$ on $\PP_R^5$ such that
  \begin{enumerate}
    \item the ${\cal G}_{a,b}$-action on each pair $x_i,y_i$ 
      is as in Lemma \ref{group action}, 
    \item such that the quadrics cutting out $\varphi(\widetilde{{\cal X}})$ 
      are $R$-linear combinations of the invariant quadrics 
      of Lemma \ref{group action}.
 \end{enumerate} 
\end{Proposition}

\prf
We consider the short exact sequence
$$
0\,\to\,\OO_{\cal X}\,\to\,\pi_\ast\OO_{\widetilde{{\cal X}}}\,\to\,\omega_{\cal X}\,\to\,0\,.
$$
Using our explicit description of $\widetilde{{\cal X}}$, we see that
${\cal G}_{a,b}$ acts on $\pi_\ast\OO_{\widetilde{{\cal X}}}$ via its regular representation
and identifies the subsheaf $\OO_{{\cal X}}$ with the ${\rm id}$-eigensheaf.

We denote the special fiber of $f$ by $X$ and the restriction
of $\overline{{\cal L}}$ by $\cal L$.
Since $h^1(X,{\cal L})=0$ by \cite[Theorem 1.5.1]{Cossec; Dolgachev},
%, semicontinuity and the fact that $\chi$
% is constant in flat families, 
global sections
of $\cal L$ extend to global sections of $\overline{{\cal L}}$.
Clearly, $\omega_{{\cal X}/R}\otimes\overline{{\cal L}}$ is an extension
of $\omega_X\otimes{\cal L}$ to ${\cal X}$
and since $h^1(X,\omega_X\otimes{\cal L})=0$, also global 
sections of $\omega_X\otimes{\cal L}$ extend to global sections of
$\omega_{\cal X}\otimes\overline{{\cal L}}$. 
In particular, $f_\ast\overline{{\cal L}}$ and 
$f_\ast(\omega_{\cal X}\otimes\overline{{\cal L}})$ are free
$R$-modules of rank $3$.

In particular, $\pi^\ast\overline{{\cal L}}$ defines a map 
$\varphi_{\overline{{\cal L}}}:\widetilde{{\cal X}}\dashrightarrow\PP^5_R$ 
that coincides with $\varphi_{\cal L}$ on the special fiber.
By Theorem \ref{three quadrics},
$\varphi_{\cal L}$ is a morphism whose image
is the complete intersection of three quadrics and so
the same is true for $\varphi_{\overline{{\cal L}}}$ by 
openness of these properties.

We take cohomology in the short exact sequence
$$
0\,\to\,\overline{{\cal L}}\,\to\,\pi_\ast\pi^\ast\overline{{\cal L}}\,\to\,
\omega_{\cal X}\otimes\overline{{\cal L}}\,\to\,0\,.
$$
Next we choose 
an $R$-basis $x_1,x_2,x_3$ of $f_\ast\overline{{\cal L}}$,
and use the ${\cal G}_{a,b}$-action on 
$f_\ast(\pi^\ast\pi_\ast{\overline{\cal L}})$ to obtain
lifts $y_1,y_2,y_3$ of an $R$-basis of 
$f_\ast(\omega_{\cal X}\otimes\overline{{\cal L}})$ to 
$f_\ast(\pi_\ast\pi^\ast\overline{{\cal L}})$ such that
${\cal G}_{a,b}$ acts on each pair
$x_i,y_i$ as in Lemma \ref{group action}.
In particular, as a ${\cal G}_{a,b}$-representation, 
$f_\ast(\pi_\ast\pi^\ast\overline{{\cal L}})$ is isomorphic 
to $\rho_{\rm reg}^{\oplus 3}$, that is,
$3$ copies of the regular representation.

Now, consider the exact sequence of ${\cal G}_{a,b}$-modules
$$
0\,\to\,\ker \mu\,\to\, {\rm Sym}^2 f_\ast\left(\pi_\ast\pi^\ast\overline{{\cal L}}\right)
\,\stackrel{\mu}{\longrightarrow}\,
f_\ast\left(\pi_\ast\pi^\ast({\overline{{\cal L}}}^{\otimes 2})\right)\,\to\,0\,.
$$
%which is exact since the occuring cohomology groups are 
%free $R$-modules and the sequence is exact modulo the
%maximal ideal of $R$.
The kernel $\ker\mu$ is easily seen to be of rank $3$.
Arguing as before, we see that 
$f_\ast(\pi_\ast\pi^\ast\overline{{\cal L}}^{\otimes 2})$ is
isomorphic to $\rho_{\rm reg}^{\oplus 9}$ as
${\cal G}_{a,b}$-representation. 
Decomposing the ${\cal G}_{a,b}$-representation on
${\rm Sym}^2 f_\ast(\pi_\ast\pi^\ast\overline{{\cal L}})$, 
we find that ${\cal G}_{a,b}$ acts trivially on $\ker\mu$.
Thus, $\varphi(\widetilde{{\cal X}})$ is cut out by three quadrics
that are ${\cal G}_{a,b}$-invariant.
In particular, these quadrics are $R$-linear combinations of the quadrics in 
Lemma \ref{group action}.
\qed\medskip

\begin{Remark}
 \label{reid bombieri mumford}
 Following an idea of Reid, Bombieri and Mumford \cite[\S3]{bm3} 
 gave the first construction of 
 all three types of Enriques surfaces in characteristic $2$.
 Our result shows that in fact {\em all} Enriques surfaces
 arise in this way - after possibly resolving Du~Val singularities
 of the quotient.
\end{Remark}

\section{Polarized Moduli and Lifting}
\label{sec:moduli}

In this section, we consider Enriques surfaces with mild singularities
together with {\em ample} Cossec--Verra polarizations.
By Proposition \ref{prop:polarization}, every
Enriques surface possesses Cossec--Verra
polarizations, but these may only be big and nef.
However, after possibly contracting nodal curves,
we obtain an Enriques surface with Du~Val singularities
together with an ample Cossec--Verra polarization.
Thus, we are led to studying pairs $(X,{\cal L})$, where $X$
is an Enriques surface with at worst Du~Val singularities
and $\cal L$ is an ample Cossec-Verra polarization.
We show that such pairs have an extremely nice
deformation theory, construct their moduli space
$\ModulPolAmple$ and prove lifting to characteristic zero.

\subsection{Enriques Surfaces with Du~Val Singularities}
Let us first slightly extend our setup and Definition \ref{def:polarization}:
A proper surface $X'$ over $k$ is an {\em Enriques surface with at worst Du~Val singularities},
if it has at worst Du~Val singularities and its minimal resolution of
singularities $\nu:X\to X'$ is a smooth Enriques surface.
Moreover, we define an invertible sheaf ${\cal L}'$ on $X'$
with at worst Du~Val singularities to be a {\em Cossec--Verra polarization}
if $\nu^\ast{\cal L}'$ on $X$ is a Cossec--Verra polarization.
Thus, if ${\cal L}$ is a Cossec--Verra polarization on 
a smooth Enriques surface $X$, then, by Theorem \ref{three quadrics},
$$
  \nu\,:\, X\,\longrightarrow\, X'\,:=\,\Proj \bigoplus_{n\geq0}H^0(X,{\cal L}^{\otimes n})
$$
is a contraction to an Enriques surface $X'$ with at worst Du~Val singularities,
and $\nu_\ast{\cal L}\,\iso\,\OO_{X'}(1)$ is an {\em ample} Cossec--Verra
polarization on $X'$.
%, which has at worst Du~Val singularities by Theorem \ref{three quadrics}.

When dealing with moduli problems, it may be necessary to consider
families of algebraic spaces rather than schemes.
Now, a smooth and proper algebraic space of dimension $2$
over a field is automatically projective, hence, a scheme, see, for example
\cite[Chapter V.4]{Knutson}.
% or: Artin, 'Algebraic Spaces' Section 4
Also, if a proper algebraic space $X'$ of dimension $2$
has rational singularities, then its minimal resolution of singularities 
will be a scheme (by what we just said), and since the contraction
to a rational singularity can always be performed in the category
of schemes \cite[Theorem (2.3)]{Artin Some numerical}, 
$X'$ is again a scheme.
In particular, over a field, the notions of Enriques surface - smooth or 
with Du~Val singularities -
in the category of schemes and algebraic spaces coincide.

Next, let us show that Enriques surfaces with Du~Val singularities 
are open in families, which is probably known to the experts.
Quite generally, in a family of varieties in characteristic zero, 
the set of fibers having only rational singularities is open  
by a result of  Elkik \cite[Th\'eor\`eme 4]{Elkik}.
His proof relies on resolution of singularities, which is not (yet) available
in positive and mixed characteristic. 
From \cite[Proposition 6.1]{Liedtke Horikawa}, we conclude

\begin{Proposition}
 \label{prop:du val open}
  Let 
  ${\cal X}\to S$ be a flat and proper family of algebraic surfaces,
  with ${\cal X}$ and $S$ schemes or algebraic spaces.
 Then, the set of points $s\in S$, such that the fiber ${\cal X}_s$ has
 rational (resp., Du~Val) singularities, is open.\qed
\end{Proposition}

\begin{Corollary}
 Under the assumptions of Proposition \ref{prop:du val open},
 the set of points $s\in S$, such that the geometric
 fiber ${\cal X}_{\overline{s}}$
 is an Enriques surface with Du~Val singularities, is open.
\end{Corollary}

\prf
If ${\cal X}_{\overline{s}}$ is an Enriques surface with Du~Val singularities,
then, by the previous result, there exists an open set $U\subseteq S$
containing $s$,
such that the fibers over points of $U$ are proper surfaces with at worst
Du~Val singularities.
By flatness, every fiber ${\cal X}_t$, $t\in U$ satisfies
$\chi(\OO_{{\cal X}_t})=\chi(\OO_{{\cal X}_s})=1$.
Since $\omega_{{\cal X}_s}^{\otimes 2}$ is trivial,
$\omega_{{\cal X}_t}$ is numerically trivial for all $t\in U$.
% have \omega in family, in one fiber it sits in Pic^\tau plus deformation...
Thus, the minimal desingularization of ${\cal X}_{\overline{t}}$
is a smooth and proper surface with $\chi=1$ and numerically
trivial canonical sheaf, that is, an Enriques surface.
% (see the introduction of \cite{bm2}).
\qed\medskip

We denote by ${\rm Aut}(X,{\cal L})$
the subgroup scheme of the automorphism group scheme ${\rm Aut}(X)$
of those automorphisms $\psi$ such that $\psi^\ast({\cal L})\iso{\cal L}$.

\begin{Proposition}
  \label{prop:automorphisms}
  Let $(X,{\cal L})$ be a Cossec--Verra polarized Enriques surface with at worst Du~Val singularities
  over  $k$.
  Then, ${\rm Aut}(X,{\cal L})$ is a group scheme of finite type over $k$.
  Moreover, if the minimal resolution of singularities of $X$ is a singular, or a classical and non-exceptional
  Enriques surface, then ${\rm Aut}(X,{\cal L})$ is finite and \'etale.
\end{Proposition}

\prf
If $\cal L$ is ample, it is well-known that ${\rm Aut}(X,{\cal L})$ is a group scheme of finite
type over $k$.
%Kollar's book and the section on the Hilbert scheme, Hom scheme etc.
If $\cal L$ is only big and nef, then $\nu:X\to X':={\rm Proj} \bigoplus_n H^0(X,{\cal L}^{\otimes n})$ 
is a contraction to an Enriques
surface with Du~Val singularities such that ${\cal L}':=\nu_\ast{\cal L}$ is ample,
and the injection ${\rm Aut}(X,{\cal L})\to{\rm Aut}(X',{\cal L}')$ gives the statement
in this case.
Let $Y$ be a minimal and smooth model of $X$.
Then, every automorphism of $X$ induces a birational self-map
of $Y$, which necessarily extends to an automorphism of $Y$ (by minimality).
Thus, we obtain an injective map ${\rm Aut}(X,{\cal L})\to {\rm Aut}(Y)$. 
If $Y$ is singular, or classical and non-exceptional, then $H^0(\Theta_{Y})=0$
(see Section \ref{sec:generalities}), that is, the identity component 
${\rm Aut}^0(Y)$ of ${\rm Aut}(Y)$ is trivial,
in which case also the identity component ${\rm Aut}^0(X,{\cal L})$
of ${\rm Aut}(X,{\cal L})$ must be trivial.
But then, being of finite type, ${\rm Aut}(X,{\cal L})$ is finite and \'etale.
\qed

\subsection{Picard scheme and effectivity}

Before studying deformations, moduli, and lifting, we have to understand 
invertible sheaves and Picard schemes in families.
In particular, the next result shows that the flatness assumption we made about 
$\Pic^\tau({\cal X}/R)$
in Section \ref{subsec:explicit} is always satisfied.

\begin{Proposition}[Ekedahl--Shepherd-Barron]
 \label{prop:picard scheme flat}
 Let $f:{\cal X}\to S$ be a flat family of Enriques surfaces with at worst Du~Val singularities,
 where ${\cal X}$ and $S$ are Noetherian algebraic spaces.
 Then,
 \begin{enumerate}
   \item $\Pic^\tau({\cal X}/S)$ is a finite flat group scheme of length $2$ over $S$.
 \end{enumerate}
 Moreover, if $f:{\cal X}\to S$ is smooth, then
 \begin{enumerate}
   \setcounter{enumi}{1}  
   \item $\Pic({\cal X}/S)$ is flat over $S$,
    \item $\Pic({\cal X}/S)/\Pic^\tau({\cal X}/S)$ is a locally constant sheaf of torsion-free
       finitely generated Abelian groups.
 \end{enumerate}
 Finally, if $f:{\cal X}\to\Spf R$ is a formal deformation of Enriques surfaces with at 
 worst Du~Val singularities, where $R$ is complete, local and Noetherian and
 with special fiber ${\cal X}_s$, then
 \begin{enumerate}
    \setcounter{enumi}{3}  
    \item the deformation is automatically projective, and hence effective, and
    \item if ${\cal L}\in\Pic({\cal X}_s)$, then ${\cal L}^{\otimes 2}$ extends to $\Spf R$.
 \end{enumerate}
\end{Proposition}

\prf
The following proof was taken from an old and unpublished draft by Ekedahl and
Shepherd-Barron, and I thank
%Since \cite{Ekedahl unpublished} is not publicly available, we reproduce 
%here a proof of a special case of  \cite[Proposition 2.2]{Ekedahl unpublished}.
% I thank 
Torsten~Ekedahl for sharing it with me.
In the meantime, this draft was largely rewritten and appeared as \cite{ESBH}, but since their article
now refers back to this article, I decided to include the following proof in order to be
self-contained and to avoid the impression of interdependencies and vicious circles
between these articles.

By the infinitesimal criterion for flatness, we 
may assume that $S=\Spec R$, where $R$ is local and Artinian, with maximal ideal
$\idealm$ and with closed point $s$, whose
residue field $k=\kappa(s)$ is algebraically closed.
Note that we have 
\begin{equation}
 \label{h2 vs h1}
  h^2(\OO_{{\cal X}_{t}})\,=\,h^1(\OO_{{\cal X}_{t}}) - \dim\Pic({\cal X}_{t})
\end{equation}
for every geometric point $t\in S$.
In particular, if ${\rm char}(k)\neq2$, then 
$h^2({\cal X}_s,\OO_{{\cal X}_s})=0$ and we are done.
Thus, we may assume that ${\rm char}(k)=2$.
%The condition 
%$h^2(\OO_{{\cal X}_{\bar{s}}})\,=\,h^1(\OO_{{\cal X}_{\bar{s}}}) - \dim\Pic^0({\cal X}_{\bar{s}})$
%implies that $H^2({\cal X}_s,W\OO_{{\cal X}_s})$ is a $W(k)$-module 
%of finite length.

Let us first show that $\Pic({\cal X}/S)$ is flat along the zero-section:
the completion $T$ of the local ring at $0\in\Pic({\cal X}/S)$ is the quotient of some
power series ring $P$ over $R$ by an ideal $\idealI$ generated by 
at most $h^2(\OO_{{\cal X}_s})$ elements, see, for example,
\cite[Lecture 27]{Mumford}.
% and Nakayama's lemma...
Thus, also the kernel $\ker(P/\idealm P\to T/\idealm T)$ is generated by
at most $h^2(\OO_{{\cal X}_s})$ elements.
Now, $\Pic({\cal X}_{\bar{s}})$ is a local complete intersection
(being a group scheme over a field), and by (\ref{h2 vs h1}) the completion
of its local ring at $0$ is a power series ring over $k$ in 
$h^1(\OO_{{\cal X}_s})$ indeterminacies modulo a complete intersection
ideal generated by $h^2(\OO_{{\cal X}_s})$ elements.
From this, we conclude that also $\idealI$ is a complete intersection ideal.
Thus, $T$ is flat over $R$.

Since $\Pic({\cal X}/S)$ is flat along the zero-section, $R^2f_\ast\widehat{\GG}_m$ 
is pro-representable by \cite[2.7.5.3]{Raynaud p-torsion}.
Now, if $f_s:{\cal X}_s\to\Spec k$ is smooth, then 
$R^2 (f_s)_\ast\widehat{\GG}_m=0$.
If ${\cal X}_s$ has at worst Du~Val singularities, then, after passing to the minimal
resolution of singularities and using the Leray spectral sequence, we find 
again $R^2 (f_s)_\ast\widehat{\GG}_m=0$.
In particular, $R^2f_\ast\widehat{\GG}_m=0$, since the restriction to the 
special fiber is zero.
To prove flatness in general, it suffices to show that any $k$-point lifts to a 
$\widetilde{S}$-point, where $\widetilde{S}\to S$ is flat, because
then, translation by this lifting induces an isomorphism of local rings.
However, the obstruction to the existence of such a lifting is an element in 
$H^2({\cal X},R^2f_\ast\widehat{\GG}_m)$,
and since $R^2f_\ast\widehat{\GG}_m=0$, this element is killed by some flat extension.

Since $R$ is Artinian and ${\cal X}_s$ is projective, $f$ is projective, and thus,
$\Pic^\tau({\cal X}/S)$ is an open sub-algebraic space.
In particular, it is also flat over $S$.
This means that $\Pic({\cal X}/S)/\Pic^\tau({\cal X}/S)$ is 
a flat group algebraic space and as it is always unramified, it is \'etale.
Since each component is proper, it is locally constant.
Also, since $\Pic^\tau$ of an Enriques surface with at worst Du~Val singularities 
is of length $2$, 
we obtain the assertion on $\Pic^\tau({\cal X}/S)$.
Moreover, if $f$ is smooth, then $\Pic/\Pic^\tau$ of an Enriques surface is 
a free Abelian group of rank $10$, and we obtain the assertion on $\Pic({\cal X}/S)/\Pic^\tau({\cal X}/S)$.

In particular, if ${\cal L}$ is an invertible sheaf on the special fiber,
then ${\cal L}^{\otimes 2}$ extends to this deformation.
Applying this to some ample ${\cal L}$,
we conclude that formal deformations are projective, and thus, effective.
\qed\medskip

\begin{Remark}
  In characteristic $\neq2$, this result is more or less trivial.
  In case $X$ is a smooth and singular (in the sense of Bombieri--Mumford) 
  Enriques surface over a perfect field
  $k$ of characteristic $2$, and $R=W(k)$,
  extension of ${\cal L}^{\otimes 4}$ 
  was established by Lang \cite[Theorem 1.4]{Lang}.
\end{Remark}

In Section \ref{sec:final}, we will also need the following, slightly technical result 
 - again, in characteristic $\neq2$ it is more or less trivial, and for singular Enriques surfaces
it was shown
in the course of the proof of \cite[Theorem 1.3]{Lang}.

\begin{Lemma}
 \label{lemma:dlog}
 Let $X$ be an Enriques surface over an algebraically closed
 field $k$ of characteristic $p>0$. 
 Then, the map
 $$
 d\log\,:\,({\rm NS}(X)/{\rm torsion})\otimes_\ZZ k\,\longrightarrow\,H^1(X,\,\Omega_X^1)
 $$
 is injective.
\end{Lemma}

\prf
Quite generally, the Chern map
$({\rm NS}(X)/{\rm tors.})\otimes_\ZZ W\to\Hcris{2}(X/W)$ is injective.
By the discusssion at the end of \cite[page 657]{Illusie},
 it induces an isomorphism
between $({\rm NS}(X)/{\rm tors.})\otimes W$
and $\Hcris{2}(X/W)/{\rm tors.}$ if $X$ is an Enriques surface.
By \cite[Corollaire II.7.3.3]{Illusie}, the slope spectral sequence of an Enriques
surface degenerates at $E_1$, and thus,
induces an isomorphism between $\Hcris{2}(X/W)/{\rm tors.}$ and
$H^1(W\Omega_X^1)/{\rm tors.}$.
Also, multiplication by $p$ induces an embedding of $(H^1(W\Omega_X^1)/{\rm tors.})\otimes_W k$
 into $H^1(\Omega_X^1)$.
% look at multiplication by p on W_n\Omega_X^i, and apply Prop. 3.4. (Illusie, page 569) 
% note that the kernel is Fil^... the kernel of the projection map to \Omega_X^i
The commutative diagram
 $$
  \xymatrixcolsep{2pc}
  \xymatrix{
      ({\rm NS}(X)/{\rm tors.})\otimes_\ZZ W \ar@{->}[r]\ar@{->}[d] & \Hcris{2}(X/W)/{\rm tors.} \,\iso\, H^1(W\Omega_X^1)/{\rm tors.} \ar@{->}[d]\\
      ({\rm NS}(X)/{\rm tors.})\otimes_\ZZ k \ar@{->}[r] & H^1(\Omega_X^1)
 }
 $$
 then shows injectivity of $d\log$, as stated.
\qed\medskip

\subsection{Polarized deformations}
\label{sec:deformation}
Next, we describe infinitesimal deformations of Cossec--Verra polarized Enriques surfaces.
Let $X$ be an Enriques surface with at worst Du~Val singularities
over $k$ and $\cal L$ be an invertible sheaf on $X$.
We define the functor
$$\begin{array}{ccccc}
  \Def_{(X,{\cal L})}&:&
\left\{  \begin{array}{l} \mbox{ local Artin algebras } \\ 
                          \mbox{ with residue field $k$ } 
         \end{array} \right\} &\to&\left(\mbox{ Sets }\right)%\\
%  && R &\mapsto&
%\left\{ \begin{array}{l} \mbox{ pairs $({\cal X},\overline{{\cal L}})$, where ${\cal X}$ is a flat}\\
%                 \mbox{ deformation of $X$ over $R$, and}\\
%                \mbox{ $\overline{\cal L}$ extends ${\cal L}$ }
%\end{array} \right\}
\end{array}
$$
that associates to each $R$ the set of pairs $({\cal X}, \overline{{\cal L}})$,
where $\cal X$ is a flat deformation of $X$ over $R$ and $\overline{{\cal L}}$
extends $\cal L$ to $\cal X$.

By Proposition \ref{prop:picard scheme flat}, $\Pic^\tau({\cal X}/R)$ is a finite and
flat group scheme of length $2$ over $R$.
By the Tate--Oort classification \cite[Theorem 2]{Tate; Oort}, it is isomorphic to
${\cal G}_{b,a}$ for some $a,b\in R$ with $ab=2$.
For a finite flat group scheme $G$ of rank $2$ over $k$ with Tate--Oort
parameters $a,b\in k$ and $ab=2$,
we denote by $\Def_G$ the functor that assigns to each local Artin algebra
$R$ with residue field $k$ the set of finite flat group schemes 
of rank $2$ over $R$ with special fiber $G$.
If $k$ is perfect, then $\Def_G$
has $W(k)[[a,b]]/(ab-2)$ as pro-representable hull. 

\begin{Theorem}
 \label{thm:unobstructed polarized deformations}
 Let $(X,{\cal L})$ be an Enriques surface with at worst Du~Val singularities together with
 an ample Cossec--Verra polarization.
 Then, the morphism of functors
 $$
 \begin{array}{ccc}
    \Def_{(X,{\cal L})} &\longrightarrow& \Def_{\Pic^\tau(X)}\,,
 \end{array}
 $$
 that assigns to a flat deformation $({\cal X},\overline{{\cal L}})/R$ its $\Pic^\tau({\cal X}/R)$,
 is formally smooth.
\end{Theorem}

\begin{Remark}
 We will see in Theorem \ref{thm:singular moduli} that this result is no longer true if
 ${\cal L}$ is not ample or when considering unpolarized deformations.
\end{Remark}

\prf
Existence of this morphism follows from Proposition \ref{prop:picard scheme flat}.

Let $R'\to R$ be a small extension, 
$({\cal X},\overline{{\cal L}})$ be a deformation of $(X,{\cal L})$
over $R$ and let ${\cal G}'$ be a finite flat group scheme extending 
${\cal G}:=\Pic^\tau({\cal X}/R)$ to $R'$.
To prove formal smoothness, we have to find an extension of 
$({\cal X}, \overline{{\cal L}})$ to $R'$ whose relative
$\Pic^\tau$ is isomorphic to ${\cal G}'$.
By Proposition \ref{three quadrics explicit},
there exists a ${\cal G}^D$-torsor
$\pi:\widetilde{{\cal X}}\to{\cal X}$ and  
$\pi^\ast(\overline{{\cal L}})$ defines an 
embedding into $\PP^5_R$, whose image is a complete 
intersection of three ${\cal G}^D$-invariant quadrics.
From the explicit description in Lemma \ref{group action}
and Proposition \ref{three quadrics explicit}, we see that
we can find a ${\cal G}'^D$-action on $\PP^5_{R'}$,
as well as 
a complete intersection of ${\cal G}'^D$-invariant quadrics
$\widetilde{{\cal X}}'$, extending $\widetilde{{\cal X}}$
together with its ${\cal G}^D$-action on $\PP^5_R$ to $R'$.

Next, we consider the map
$\Psi:=\Psi^{{\cal G}'^D}: \PP^5_{R'} \to \PP_{R'}^{11}$
%$$
%  \Psi\,:\,\PP^5_{R'} \,\to\,\PP_{R'}^{11}
%$$
given by the ${\cal G}'^D$-invariant quadrics of 
Lemma \ref{group action}, which is easily seen to be a morphism.
By the same lemma, the ${\cal G}'^D$-invariants 
of even degree in $R'[x_1,y_1,...,x_3,y_3]$ are generated
by these $12$ invariant quadrics.
%\marginpar{if Lemma fails use semi-continuity!}
Thus, we can identify $\Psi:\PP^5_{R'}\to \Psi(\PP^5_{R'})$ with
the quotient morphism by ${\cal G}'^D$.
In particular, 
${\cal X}':=\Psi(\widetilde{{\cal X}}')\iso{\widetilde{\cal X}}'/{\cal G}'^D$ 
is flat over $R'$, extending $\widetilde{{\cal X}}/{\cal G}^D\iso{\cal X}$
to $R'$.

Finally, the ${\cal G}'^D$-action defines a descent data on
$\OO_{\PP^5_{R'}}(1)|_{\widetilde{{\cal X}}'}$.
Thus, by finite flat descent, there exists an invertible sheaf
${\cal L}'$ on ${\cal X}'$ that extends ${\cal L}$.
\qed\medskip

\subsection{Lifting to characteristic zero}

As an application of the results established so far,
we show that Enriques surfaces lift to characteristic zero.
In the case of polarized Enriques surfaces,
the ramification is completely controlled by $\Pic^\tau$:

\begin{Theorem}
 \label{thm:nice lifting}
 Let $X$ be an Enriques surface with at worst Du~Val singularities.
 Assume that $X$ admits an ample Cossec--Verra polarization. 
 \begin{enumerate}
  \item If $X$ is not supersingular, then it lifts over the Witt ring $W(k)$,
  \item if $X$ is supersingular, then it lifts over $W(k)[\sqrt{2}]$, but not over $W(k)$.
 \end{enumerate}
 Here, all lifts are projective, and in particular, algebraizable.
\end{Theorem}

\prf
We fix $a,b\in k$ with $ab=2$ such that $\Pic^\tau(X)\iso G_{b,a}$.
%(Tate--Oort parameters).
If $X$ is not supersingular then $a\neq0$ or $b\neq0$ and we can find
$a',b'\in W(k)$ with $a'b'=2$ lifting $a,b$, respectively, and then
${\cal G}_{b',a'}$ defines a lift of $G_{b,a}$ to $W(k)$.
Lifting of $X$ over $W(k)$ then follows by running through the proof
of Theorem \ref{thm:unobstructed polarized deformations} with
$R=k$ and $R'=W(k)$.
Of course, $W(k)\to k$ is not a small extension but the proof also
works in this case.
(Alternatively, Theorem \ref{thm:unobstructed polarized deformations}
provides us with a formal lifting of $X$ over $\Spf W(k)$,
which is algebraizable by Proposition \ref{prop:picard scheme flat}.)

If $X$ is supersingular, then $a=b=0$ and ${\cal G}_{\pi,\pi}$ with
$\pi:=\sqrt{2}$ lifts $G_{0,0}\iso\alpha_2$ to $W(k)[\sqrt{2}]$.
The same arguments as before show that $X$ lifts over $W(k)[\sqrt{2}]$.
On the other hand, by Proposition \ref{prop:picard scheme flat},
if ${\cal X}$ were a lifting of $X$ over $W(k)$
then $\Pic^\tau({\cal X}/R)$ would be a finite flat group scheme
of length $2$ over $W(k)$ with special fiber $\alpha_2$.
However, this contradicts the fact that $\alpha_p$ does not admit
liftings over $W(k)$, see also \cite[Theorem 2]{Tate; Oort}.
\qed\medskip

Since a smooth Enriques surface admits an ample
Cossec--Verra polarization after possibly contracting
$(-2)$-curves to Du~Val singularities, the previous
theorem combined with Artin's results \cite{Artin Brieskorn}
on simultaneous resolutions of surface singularities in families
implies the following lifting result for smooth Enriques surfaces:

\begin{Theorem}
 \label{thm:general lifting}
 Let $X$ be an Enriques surface % in positive characteristic $p$.
 over an algebraically closed field $k$
 of positive characteristic $p$.
 \begin{enumerate}
  \item If $p\neq2$, then $X$ lifts over $W(k)$.
  \item If $X$ is a singular, or a classical and non-exceptional
    Enriques surface, then it lifts over $W(k)$.
  \item If $X$ is supersingular and admits an ample Cossec--Verra 
    polarization then it lifts over $W(k)[\sqrt{2}]$, but not over $W(k)$.
  \item In the remaining cases, $X$ lifts (as special fiber of an algebraic space)
    over a possibly ramified extension 
    of $W(k)$.
 \end{enumerate}
 Here, all lifts are projective, and in particular, algebraizable.
\end{Theorem}

\prf
In the first two cases, we have $h^2(\Theta_X)=0$, that is, there
exists a formal lifting over $W(k)$, whose algebraization follows
from Proposition \ref{prop:picard scheme flat}.
Lifting over $W(k)[\sqrt{2}]$ and non-lifting over $W(k)$
for supersingular surfaces admitting
ample Cossec--Verra polarizations has been established in 
Theorem \ref{thm:nice lifting}.
In the remaining cases, we find a birational model of $X$ with at
worst Du~Val singularities that admits an ample Cossec--Verra
polarization and that lifts according to Theorem \ref{thm:nice lifting}.
By \cite[Theorem 3]{Artin Brieskorn}, there exists
an algebraic space, smooth over a possibly ramified extension of $W(k)$,
with special fiber $X$.
Algebraizability follows again from Proposition \ref{prop:picard scheme flat}.
\qed

\begin{Remark}
 The first two results were already known, see \cite{Lang} and \cite{Ekedahl}.
 It would be interesting to know the ramification needed
 to lift exceptional, as well as supersingular Enriques surfaces.
 Some bounds on the ramification are given in \cite[Corollary 5.7]{ESBH}.
 By Proposition \ref{prop:polarization},
 Enriques surfaces need not possess ample Cossec--Verra polarizations. 
 In fact, we will see in Theorem \ref{thm:singular moduli} and Remark \ref{resolution bad} 
 that the moduli space at exceptional Enriques surfaces may be more singular
 than expected.
\end{Remark}

\subsection{Moduli spaces for Cossec--Verra polarized Enriques surfaces}
For a Noetherian base $S$, we consider proper and flat morphisms of
algebraic spaces ${\cal X}\to S$ together with ${\cal L}\in\Pic({\cal X}/S)$, flat
over $S$,
such that every geometric fiber is an Enriques surface with at worst
Du~Val singularities, and such that ${\cal L}$ restricts to an ample
Cossec--Verra polarization on every geometric fiber.
We denote the set of all such pairs $({\cal X}/S)$ by
$\ModulPolAmple$, that is, consider the moduli problem
$$
\begin{array}{cl}
& \mbox{ \quad $S$-valued points} \\
\hline

\ModulPolAmple & \mbox{morphisms } ({\cal X},{\cal L})\to S\mbox { of algebraic spaces, whose geometric fibers}\\
&\mbox{are Enriques surfaces with at worst Du~Val singularities, and where }\\
&\mbox{${\cal L}\in\Pic(X/S)$ is an invertible sheaf that restricts to an ample}\\
&\mbox{Cossec--Verra polarization on each geometric fiber}
\end{array}
$$
Next, we show that this set carries the structure of an
Artin stack, and describe its geometry:

\begin{Theorem}
  \label{thm:artin stack}
  $\ModulPolAmple$ carries the structure of a quasi-separated Artin stack of finite type 
  over $\Spec \ZZ$, and
  $\ModulPolAmple\otimes_\ZZ\ZZ[\frac{1}{2}]$ is even Deligne--Mumford. 
  For a field $k$ of characteristic $p$
  \begin{enumerate}
   \item If $p\neq2$, then $\ModulPolAmple\otimes_\ZZ k$ is irreducible, unirational, $10$-dimensional and 
     smooth over $k$.
   \item If $p=2$, then $\ModulPolAmple\otimes_\ZZ k$ consists of two irreducible, unirational, smooth 
     and $10$-dimensional components 
     $$
       \MMod^{\mu_2} \mbox{ \quad and \quad }\MMod^{\ZZ/2\ZZ}\,.
     $$
     Moreover,
     \begin{enumerate}
     \item[-] they intersect transversally along an irreducible, unirational, smooth and $9$-dimensional 
        closed substack $\MMod^{\alpha_2}$,
     \item[-] $\MMod^{\alpha_2}$ parametrizes supersingular surfaces,
     \item[-] $\MMod^G\,-\,\ModulPolAmple^{\alpha_2}$ parametrizes
           singular surfaces ($G=\mu_2$) and 
           classical surfaces ($G=\ZZ/2\ZZ$), respectively,
     \item[-] for all $G$, $\MMod^G$ 
           contains an open and dense substack, whose geometric points correspond to
           smooth surfaces.
     \end{enumerate}
  \end{enumerate}
\end{Theorem}

\prf
First of all, since we consider Cossec--Verra polarized families, every formal deformation 
is effective.
(Alternatively, we could also use Proposition \ref{prop:picard scheme flat}.)
Under our assumptions, it is standard that $\ModulPolAmple$ 
can be given the structure of a quasi-separated Artin stack of finite type over $\ZZ$, 
see \cite[Example (5.5)]{Artin versal deformations} for a 
sketch (the role of the canonical polarization in Artin's example is replaced by
Cossec--Verra polarizations in our setup, and the results on automorphisms
are provided by Proposition \ref{prop:automorphisms}), or \cite{rizov} for a detailed and 
elaborate discussion.
Since the automorphism group schemes are finite and \'etale in characteristic $\neq2$,
$\ModulPolAmple\otimes_\ZZ\ZZ[\frac{1}{2}]$ is even Deligne--Mumford.

Now, we base-change to a field $k$ of characteristic $p\geq0$.
We only discuss characteristic $2$, since the analysis for $p\neq2$ is analogous
to the case of singular Enriques surfaces in characteristic $2$.

Let $G$ be a finite flat group scheme of length $2$ over $k$
and consider the $G^D$-action on $\PP_k^5$ defined in 
Proposition \ref{three quadrics explicit}.
As in the proof of Theorem \ref{thm:unobstructed polarized deformations},
we define $\Psi^G:\PP_k^5\to\PP_k^{11}$ to be the morphism defined by 
the $G^D$-invariant quadrics of Lemma \ref{group action}.
As explained in \cite[\S 3]{bm3}, the inverse image of a geometric
generic hyperplane
section yields the canonical double cover of an Enriques 
surface together with a $G^D$-action.
Such surfaces are overparametrized by an open and dense subset $U_G$
of the Gra\ss mannian ${\rm Gr}(3,12)$ of linear $3$-dimensional
subspaces of a $12$-dimensional vector space.
By Proposition \ref{three quadrics explicit},
all Enriques surfaces with at worst
Du~Val singularities, equipped with an ample Cossec--Verra polarization,
and whose canonical double cover is a $G^D$-torsor, arise this way.
If we denote by ${\mathscr N}^G$ the substack of surfaces
with $\Pic^\tau\iso G$, then we 
have just shown that $U_G$ maps dominantly onto 
${\mathscr N}^G$, 
showing its irreducibility, as well as its unirationality.

Now, if ${\cal X}\to S$ is a family of Enriques surfaces with at worst
Du~Val singularities, then the set of points such that ${\cal X}_s$ 
is a classical Enriques surface is open, since the property
$h^1(\OO_X)=0$ is open by semi-continuity.
Also, the set of points such that ${\cal X}_s$ is singular is open:
by \cite{Artin Brieskorn} there exists a surjective map
$S'\to S$ and an algebraic space ${\cal X}'\to S'$ that simultaneously
resolves the singularities of ${\cal X}\to S$.
But then, the property of being a singular Enriques surface, that is,
satisfying $h^0(\Omega_X)=0$, is open on $S'$ by semi-continuity.
Now, for each $s\in S$, the map on Henselizations
$\OO_{S,s}^h\to\OO_{S',s'}^h$ is finite.
We conclude that being a singular Enriques surface is stable under
generization also on $S$, proving openness.
Similar arguments show that the set of points 
such that ${\cal X}_s$ is supersingular, is closed. 
We conclude that ${\mathscr N}^G$ for $G=\mu_2$ and for
$G=\ZZ/2\ZZ$ belong to different components of $\ModulPolAmple$.
We denote these components by $\MMod^G$ and remark that
they contain ${\mathscr N}^G$ as open and dense substacks.

Next, let $(X,{\cal L})\to\Spec\overline{K}$ be a geometric point of 
${\mathscr N}^{\alpha_2}$.
The finite flat group scheme ${\cal G}_{0,t}$ of length $2$ over $\overline{K}[[t]]$ 
(in the Tate--Oort notation from \cite{Tate; Oort}) has special fiber
$\alpha_2$ and generic fiber $\mu_2$.
By Theorem \ref{thm:unobstructed polarized deformations},
there exists a family ${\cal X}\to\Spec \overline{K}[[t]]$ of 
Cossec--Verra polarized Enriques surfaces with at worst Du~Val singularities
with special fiber $(X,{\cal L})$ and $\Pic^\tau({\cal X}/\overline{K}[[t]])\iso{\cal G}_{0,t}$.
In particular, the geometric generic fiber
of this family is a singular Enriques surface.
This shows that ${\mathscr N}^{\alpha_2}$ is a closed substack of 
$\MMod^{\mu_2}$.
Using ${\cal G}_{t,0}$ instead of ${\cal G}_{0,t}$ in the previous
discussion, we conclude that ${\mathscr N}^{\alpha_2}$ 
is also a closed substack of $\MMod^{\ZZ/2\ZZ}$.
Thus, $\MMod^{\mu_2}$ and $\MMod^{\ZZ/2\ZZ}$
intersect along ${\mathscr N}^{\alpha_2}$.
To show that the intersection is transversal in ${\mathscr N}^{\alpha_2}$,
we restrict the deformation functor $\Def_{\alpha_2}$ of $\alpha_2=G_{0,0}$ 
to the subcategory of local Artin $k$-algebras.
This has $k[[x,y]]/(xy)$ pro-representable hull, see the discussion
at the beginning of Section \ref{sec:deformation}.
But then,
the statements about smoothness and transversal intersections
follow from Theorem \ref{thm:unobstructed polarized deformations}.

As shown in \cite[\S3]{bm3},
generic geometric hyperplane sections of 
$\Psi^G(\PP_{k}^5)\subseteq\PP_{{k}}^{11}$
yield smooth singular ($G=\ZZ/2\ZZ$), classical ($G=\mu_2$)
and supersingular ($G=\alpha_2$) Enriques surfaces, respectively.
Clearly, $\OO_{\PP^{11}}(1)$ restricts to an ample Cossec--Verra
polarization on these surfaces.
Since ampleness and smoothness are open properties,
there exist open and dense
substacks of $\MMod^G$ for $G=\mu_2,\ZZ/2\ZZ$, and $\alpha_2$,
respectively, whose geometric points correspond to smooth surfaces.

We postpone the computation of the dimension of the three components
to the proofs of Theorem \ref{thm:nice moduli} and the proof of Theorem \ref{thm:singular moduli} 
below:
namely, we shall prove there 
that the just-established open and dense substacks of $\MMod^G$
parametrizing smooth surfaces 
are isomorphic to open substacks of the yet to be defined
stack $\ModulPolSmooth$. 
Thus, it will suffice to compute the dimension at such
points of $\ModulPolSmooth$, which is equal to $10$.
%From this, the dimension and the local description
%of $\ModulPolAmple^{\alpha_2}$ follows from 
%Theorem \ref{thm:unobstructed polarized deformations} and
%the description of the pro-representable hull of
%$\Def_{\alpha_2}$.
\qed

\begin{Remark}
 Over the complex numbers, Casnati \cite{Casnati}
 considered degree $4$ polarized Enriques surfaces, and
 showed that the corresponding moduli space is rational.
 Clearly, not every such polarization is Cossec--Verra, only generic ones.
 In view of Casnati's result,
 it would be interesting to know whether the 
 components of $\ModulPolAmple\otimes_\ZZ k$ are rational for every field $k$.
\end{Remark}

%As in Theorem \ref{thm:unobstructed polarized deformations},
Given a family ${\cal X}\to S$ of $\ModulPolAmple$, we set
$G:=\Pic^\tau({\cal X}/S)$.
Then, the augmentation ideal ${\cal I}:=\ker(\OO_G\to\OO_S)$ is an invertible
sheaf on $S$.
Zariski locally, we can trivialize ${\cal I}$ and every trivialization yields
two elements $a,b\in H^0(S,\OO_S)$ with $ab=2$,
see the introduction of \cite{Tate; Oort}.
In terms of the Tate--Oort classification \cite[Theorem 2]{Tate; Oort},
this means that we have fixed, Zariski locally on $S$, an isomorphism 
$G\iso {\cal G}_{b,a}$.
Changing the trivialization of $\cal I$,
has the effect $a\mapsto sa$, $b\mapsto s^{-1}b$
for some $s\in H^0(S,\OO_S^\times)$ on Tate--Oort parameters, defining a $\GG_m$-action.
Thus, we obtain a morphism of stacks
$$
\ModulPolAmple\,\longrightarrow\, \left[ \left(\Spec\ZZ[a,b]/(ab-2)\right)\, /\, \GG_m \right].
$$
This said, we leave the following straight-forward generalization of 
Theorem \ref{thm:artin stack} to the reader.

\begin{Theorem}
 \label{thm:not needed}
 $\ModulPolAmple$ is a
 quasi-separated Artin stack of finite type,
 smooth, and of relative dimension $10$ 
 over $[ (\Spec\ZZ[x,y]/(xy-2))/\GG_m]$.
 Moreover, $\ModulPolAmple\otimes_\ZZ\ZZ[\frac{1}{2}]$
 is a quasi-separated Deligne--Mumford stack of finite type,
 smooth and of relative dimension $10$ over $\Spec\ZZ[\frac{1}{2}]$.
 \qed
\end{Theorem}

\section{Moduli of Smooth Enriques Surfaces}
\label{sec:final}

In the previous section, we constructed and described 
$\ModulPolAmple$, the moduli space of Enriques surfaces with
at worst Du~Val singularities together with an ample Cossec--Verra
polarization.
In this section, we consider the following related moduli spaces
that parametrize smooth Enriques surfaces:
$$
\begin{array}{cl}
& \mbox{ \quad $S$-valued points} \\
\hline

\Modul & \mbox{morphisms }{\cal X}\to S\mbox{ of algebraic spaces, whose}\\
&\mbox{geometric fibers are smooth Enriques surfaces}\\
\ModulPolSmooth & \mbox{morphisms } ({\cal X},{\cal L})\to S\mbox { of algebraic spaces, whose}\\
&\mbox{geometric fibers are smooth Enriques surfaces, and }\\
&\mbox{where $\cal L$ is an invertible sheaf that restricts to a }\\
&\mbox{Cossec--Verra polarization on each geometric fiber}
\end{array}
$$
%Let us consider the following sets:
%\begin{enumerate}
% \item $\Modul$, whose elements are smooth Enriques surfaces.
% \item $\ModulPolAmple$, whose elements are pairs $(X,{\cal L})$, where
%   \begin{enumerate}
%       \item[-] $X$ is an Enriques surface with at worst Du~Val singularities, and
%       \item[-] $\cal L$ is an ample Cossec--Verra polarization.
%    \end{enumerate}
% \item $\ModulPolSmooth$, whose elements are pairs $(X,{\cal L})$, where
%    \begin{enumerate}
%       \item[-] $X$ is a smooth Enriques surface, and
%       \item[-] $\cal L$ is a Cossec--Verra polarization.
%    \end{enumerate}
% \item $\ModulGrp$, whose elements are finite flat group schemes of length $2$.
%\end{enumerate}

\subsection{Basic properties and functors}
Clearly, both moduli spaces are stacks.
Moreover, they are related by two functors:
$$
 \xymatrixcolsep{3pc}
 \xymatrix{
      & \ar[dl]_{\Phi_{\rm cont}} \ar[dr]^{\Phi_{\rm forget}} \ModulPolSmooth\\
      \ModulPolAmple & & \Modul }
$$
First, we define $\Phi_{\rm cont}$ to be the following contraction functor:
$$\begin{array}{ccc}
  \ModulPolSmooth & \to & \ModulPolAmple\\
     \left( f:({\cal X},{\cal L})\to S \right) &\mapsto& 
      \left( ({\cal X}':=\Proj \bigoplus_{n\geq0} f_\ast\left( {\cal L}^{\otimes n}) \right),\,\OO_{{\cal X}'}(1)) 
      \to S \right)
  \end{array}$$
Now, a morphism $S\to \ModulPolAmple$ corresponds 
to a family ${\cal Y}\to S$ of Enriques surfaces with at worst
Du~Val singularities together with an invertible sheaf $\cal F$ that
restricts to an ample Cossec--Verra polarization on every geometric fiber.
Thus,  the fiber of $\Phi_{\rm cont}$ over $S$ is given by the set-valued functor
$$\begin{array}{ccccc}
  \Phi_{\rm cont}^{-1}(S) &:& S'/S & \mapsto &
  \left\{
  \begin{array}{l} 
     \mbox{ simultaneous resolutions of singularities of } \\
     {\cal Y}\times_S S'\to S' \mbox{ together with the pullback of } 
     {\cal F}
  \end{array}
  \right\}
\end{array}
$$
By a result of Artin \cite[Theorem 1]{Artin Brieskorn}, this functor is representable by a 
locally quasi-separated and quasi-finite algebraic space over $S$.
Thus, 

\begin{Proposition}
  \label{prop:contraction functor}
  The functor $\Phi_{\rm cont}$ is representable and locally quasi-separated.
  It induces a bijection on geometric points. \qed
\end{Proposition}

Combined with Theorem \ref{thm:artin stack}, this gives us more structure:

\begin{Corollary}
% \label{thm:quasi-separated} 
 The moduli space $\ModulPolSmooth$ is a quasi-separated Artin stack of finite type over 
 $\Spec \ZZ$.
 Moreover, $\ModulPolSmooth\otimes_\ZZ\ZZ[\frac{1}{2}]$ is a quasi-separated Deligne--Mumford stack.
\end{Corollary}

\begin{Remark}
 \label{rem:not quasi-separated}
 Already over the complex numbers, the automorphism group of a generic Enriques
 surface is infinite and discrete (see, \cite{Dolgachev Automorphisms}, for example),  
 and so, although $\Modul$ is a stack, its diagonal will not be quasi-compact. 
 However, having a quasi-compact diagonal (quasi-separatedness)
 is usually built into the theory of algebraic stacks from the very beginning
 (as in \cite[Definition 4.1]{Laumon}, for example).
 Therefore, we will be careful with statements about the ``geometry'' of 
 $\Modul$ in the sequel.
 \end{Remark}

Next, we consider the forgetful functor 
$$\begin{array}{ccccc}
    \blowup \Phi_{\rm forget} &:& \ModulPolSmooth & \to & \Modul\\
     & & \left( ({\cal X},{\cal L})\to S \right) &\mapsto& ( {\cal X}\to S )
\end{array}$$
It has the following basic properties:

\begin{Proposition}
 \label{prop:forgetful functor}
 The functor $\Phi_{\rm forget}$ is representable, separated and
 locally finite and flat. 
 It is smooth at a geometric point $(X,{\cal L})$ of $\ModulPolSmooth$  
 if and only if $X$ is classical, that is, if and only if
 $\Pic^0(X)$ is reduced.
\end{Proposition}

\prf
%By \cite[\numero 232, Theorem 3.1]{fga},
For a family ${\cal X}\to S$ of Enriques surfaces over a Noetherian base,
$\Pic({\cal X}/S)$ is representable
by a separated algebraic space that is locally of finite type over $S$
by Proposition \ref{prop:picard scheme flat}.
It is not difficult to see that
Cossec--Verra polarizations lie discrete in the N\'eron--Severi lattice,
and are open in families.
%The condition ${\cal L}^2=4$ is closed and such invertible sheaves
%lie discrete in the N\'eron--Severi lattice.
%From Theorem \ref{sd1} and Theorem \ref{sd3} we see that the
%$\Phi\geq2$ is characterized by the property that $\varphi_{{\cal L}}$
%is birational.
%Since this is an open property, so is the condition $\Phi\geq2$.
%On the other hand,  
%invertible sheaves with self-intersection number $4$ 
%satisfy $\Phi\leq2$ by \cite[Lemma 3.6.1]{Cossec; Dolgachev}, and thus
%$\Phi=2$ is an open property for such invertible sheaves.
Thus, $\Phi_{\rm forget}^{-1}(S)$ is locally represented 
by $\Pic^\tau({\cal X}/S)$, proving local finiteness and 
flatness of $\Phi_{\rm forget}$.
In particular, $\Phi_{\rm forget}$ is smooth at a geometric point
$(X,{\cal L})$ of $\ModulPolSmooth$
if and only if $\Pic^0(X)$ is reduced.
\qed\medskip

\subsection{Local and Global Structure of the Moduli Spaces}
Now, we study the deformation theory of polarized Enriques surfaces:
quite generally, an invertible sheaf $\cal L$ on a smooth variety $X$ determines
via $d\log$ a class in $H^1(X,\Omega_X)$, the {\em Chern class} of $\cal L$.
This can be interpreted as an extension class in 
${\rm Ext}^1(\Theta_X,\OO_X)$, the {\em Atiyah extension}
of $\cal L$
$$
0\,\to\,\OO_X\,\to\,{\cal A}_{\cal L}\,\to\,\Theta_X\,\to\,0\,.
$$
The  groups $H^i(X,\Theta_X)$ provide a tangent-obstruction
theory for deformations of $X$, and similarly, 
the groups  $H^i(X,{\cal A}_{\cal L})$ provide a tangent-obstruction theory 
for deforming the pair $(X,{\cal L})$.
In particular, the differential of $\Phi_{\rm forget}$ can be computed from 
the cohomology sequence
$$
...\,\to\,H^1(\OO_X)\,\to\,H^1({\cal A}_{\cal L})
\,\stackrel{d\Phi_{\rm forget}}{\longrightarrow}\,H^1(\Theta_X)\,\to\,
H^2(\OO_X)\,\to\,...
$$
The following result is crucial for the local structure of
our moduli spaces:

\begin{Lemma}
 If $(X,{\cal L})$ is an Enriques surface with Cossec--Verra polarization ${\cal L}$,
 then we have
$$
   h^0({\cal A}_{\cal L})=h^0(\Theta_X)+1,\mbox{ \quad }
   h^1({\cal A}_{\cal L})=h^1(\Theta_X),\mbox{ \quad and \quad }
   h^2({\cal A}_{\cal L})=h^2(\Theta_X).
$$
 The values $h^i(\Theta_X)$ are well-known (see Section \ref{sec:generalities}). 
\end{Lemma}

\prf
If $X$ is classical, that is, $h^i(\OO_X)=0$ for $i\geq1$, then this is trivial.

If $X$ is non-classical, then $\Pic^0(X)$ is non-reduced.
Now, $\Phi_{\rm forget}$ is locally a torsor under $\Pic^0(X)$, which shows that
the differential $d\Phi_{\rm forget}$ is not injective,
and thus, $H^1(\OO_X)$ injects into $H^1({\cal A}_{\cal L})$.
We are done once we show that $d\Phi_{\rm forget}$ is not surjective.
Now, since ${\cal L}^2=4$, it follows that ${\cal L}$ is not $2$-divisible in $\Pic(X)$.
In particular, the class $d\log ({\cal L})\in H^1(\Omega_X^1)$
is non-zero by Lemma \ref{lemma:dlog}.
Since $\omega_X\iso\OO_X$, the pairing 
$$
H^1(\Omega_X^1)\times H^1(\Theta_X)\,\stackrel{\cup}{\longrightarrow}
H^2(\Omega_X^1\otimes \Theta_X)\,\to\,H^2(\OO_X)
$$
coincides with Serre duality, which is perfect.
% Sernesi, page 150 below
In particular, there exists a deformation ${\cal X}$ of $X$ to
$k[\epsilon]/\epsilon^2$ with a Kodaira--Spencer class $\xi$
that has the property that $d\log ({\cal L})\cdot \xi\neq0$ under the
above pairing.
Thus, ${\cal L}$ does not extend to ${\cal X}$ showing that
$d\Phi_{\rm forget}$ is not surjective.
\qed\medskip

Using the previously defined functors $\Phi_{\rm cont}$ and $\Phi_{\rm forget}$, 
we now compare our three moduli spaces $\Modul$, $\ModulPolSmooth$, and $\ModulPolAmple$.
This extends the results of Theorem \ref{thm:artin stack}.
Let us do the case of characteristic $\neq2$ first:

\begin{Theorem}
 \label{thm:nice moduli}
 The stacks
 $\ModulPolSmooth\otimes_\ZZ \ZZ[\frac{1}{2}]$
 and $\ModulPolAmple\otimes_\ZZ \ZZ[\frac{1}{2}]$ are smooth, irreducible, unirational and 
 $10$-dimensional over $\Spec\ZZ[\frac{1}{2}]$.
\end{Theorem}

\prf
We have shown this for $\ModulPolAmple\otimes_\ZZ \ZZ[\frac{1}{2}]$ in Theorem \ref{thm:artin stack}, except
for the dimension.
The assertions on smoothness and dimension of $\ModulPolSmooth$
follow from deformation theory 
and $h^1({\cal A}_{\cal L})=h^1(\Theta_X)=10$
and $h^2({\cal A}_{\cal L})=h^2(\Theta_X)=0$,
whenever $X$ is a geometric point with residue characteristic $\neq2$.
For an algebraically closed field of characteristic $\neq2$, there exists an open and dense substack, over which 
$\ModulPolSmooth\otimes_\ZZ k$ and $\ModulPolAmple\otimes_\ZZ k$ are isomorphic via
$\Phi_{\rm cont}$.
This shows irreducibility and unirationality of $\ModulPolSmooth\otimes_\ZZ\ZZ[\frac{1}{2}]$, as well
as the dimension of $\ModulPolAmple\otimes_\ZZ \ZZ[\frac{1}{2}]$.
\qed\medskip

\begin{Remark}
 Despite the discussion in Remark \ref{rem:not quasi-separated}, one could use $\Phi_{\rm forget}$
 and Proposition \ref{prop:forgetful functor} to argue that $\Modul\otimes_\ZZ\ZZ[\frac{1}{2}]$ is 
 smooth, unirational, and $10$-dimensional.
 Over the complex numbers, an analytic moduli space of unpolarized Enriques
 surfaces can be constructed via a period map, a period domain and the
 Torelli theorem, see \cite[Section VIII.20]{bhpv} for references and details.
 Using analytic methods, Kond\={o} \cite{Kondo} has shown that 
 this moduli space is rational.
% Together with our lifting results and 
% Matsusaka's Theorem \cite{Matsusaka} this implies that all moduli spaces 
% considered in Theorem \ref{moduli theorem} are birationally ruled.
 It would be interesting to extend this to positive characteristic.
 Finally, we note that Torelli theorems and period maps are not available in positive
 or mixed characteristic, although first steps are taken in 
 \cite{ESBH}.
\end{Remark}

We end our article by dealing with the  structure of our three moduli spaces
in characteristic $2$.
The proof of the following result is completely analogous to that of Theorem \ref{thm:nice moduli},
which is why we leave it to the reader:

\begin{Theorem}
 \label{thm:general moduli}
  Let $\MMod$ be equal to
  $\ModulPolAmple\otimes_\ZZ \FF_2$ or $\ModulPolSmooth\otimes_\ZZ \FF_2$.
  Then, $\MMod$ consists of two $10$-dimensional, irreducible and unirational components
  $$
          \MMod^{\mu_2}\mbox{ \quad and \quad }\MMod^{\ZZ/2\ZZ}.
  $$
  Moreover,
  \begin{enumerate}
   \item[-] they intersect along a closed substack $\MMod^{\alpha_2}$, 
        which is $9$-dimensional, irreducible and unirational,
   \item[-] $\MMod^{\alpha_2}$ parametrizes supersingular Enriques surfaces,
   \item[-] $\MMod^{G}-\MMod^{\alpha_2}$ parametrizes singular ($G=\mu_2$)
        and classical ($G=\ZZ/2\ZZ$) Enriques surfaces, respectively.\qed
  \end{enumerate}
 \end{Theorem}
 
 %Again, a similar statement can be made about $\Modul\otimes_\ZZ\FF_2$, taken with a
 %grain of salt, and bearing
 %Remark \ref{rem:not quasi-separated} in mind.
 %
The local structure of the moduli spaces is given by the following result:

 \begin{Theorem} 
 \label{thm:singular moduli}
 Let $(X,{\cal L})$ be a geometric point of $\ModulPolSmooth\otimes_\ZZ\FF_2$.
 Then,
 \begin{enumerate}
  \item[-] If $X$ is singular, or classical and not exceptional, then
    both moduli stacks
    are smooth at $(X,{\cal L})$, and $\Phi_{\rm cont}(X,{\cal L})$, 
    respectively.
  \item[-] If $X$ is a classical and exceptional Enriques surface, then
    \begin{enumerate}
     \item[-] $\ModulPolAmple$ is smooth at $\Phi_{\rm cont}(X,{\cal L})$, whereas
     \item[-] $\ModulPolSmooth$ is not smooth at 
       $(X,{\cal L})$.
 \end{enumerate}
  \item[-] If $X$ is supersingular, then the intersection of $\MMod^{\mu_2}$ and 
      $\MMod^{\ZZ/2\ZZ}$ in $X$ is transversal 
      \begin{enumerate}
        \item[-] at $\Phi_{\rm cont}(X,{\cal L})\,\in\,\ModulPolAmple$, and
        \item[-] at $(X,{\cal L})\in\ModulPolSmooth$ if $\cal L$ is ample.
      \end{enumerate}
      The latter two conditions hold along an open and dense substack of $\MMod^{\alpha_2}$.
   \end{enumerate}
\end{Theorem}

\prf
For $\ModulPolAmple\otimes_\ZZ \FF_2$, we have shown all assertions in
Theorem \ref{thm:artin stack}.
Moreover, if $X$ is singular, or classical and not exceptional, then
$h^2({\cal A}_{\cal L})=0$, and we get  smoothness 
around the corresponding points of $\ModulPolAmple\otimes_\ZZ\FF_2$.

If $X$ is a classical and exceptional Enriques surface, then 
$h^1({\cal A}_{\cal L})=12$.
Thus, if deformations of $(X,{\cal L})$ were unobstructed, (that is, the obstruction class
in the one-dimensional space $H^2({\cal A}_{\cal L})$ is in fact zero,) then
$\ModulPolSmooth$ would be 
$h^1({\cal A}_{\cal L})-h^0({\cal A}_{\cal L})=11$-dimensional at $X$.
However, $\ModulPolSmooth\otimes_\ZZ\FF_2$ is $10$-dimensional
at this point, and we conclude that  $\ModulPolSmooth$ cannot be smooth 
at  $(X,{\cal L})$.

Now, let $X$ be a supersingular Enriques surface.
Transversality of the intersection of the two components
at $\Phi_{\rm cont}(X,{\cal L})\in\ModulPolAmple$
has been established in Theorem \ref{thm:artin stack}.
Using $\Phi_{\rm cont}$, we conclude
that also $\ModulPolSmooth\otimes_\ZZ\FF_2$ consists of two components
intersecting along the supersingular locus.
Now, in case ${\cal L}$ is ample, then $\ModulPolSmooth$ and
$\ModulPolAmple$ are locally isomorphic near $(X,{\cal L})$,
%(via $\Phi_{\rm cont}$),
and thus, the intersection is transversal.
\qed\medskip

\begin{Remark}
A picture similar to that of $\ModulPolSmooth\otimes_\ZZ\FF_2$ also emerges for
$\Modul\otimes_\ZZ\FF_2$, taken with a grain of salt, and bearing
Remark \ref{rem:not quasi-separated} in mind.
\end{Remark}

We end this article by a couple of remarks concerning the rather unexpected and 
surprising phenomenon of non-smoothness 
of $\ModulPolSmooth$ (and $\Modul$) at points corresponding to exceptional Enriques
surfaces in characteristic $2$:
 
\begin{Remark}
 \label{resolution bad}
 The hull of the local deformation functor of an exceptional Enriques surface was computed in
 \cite[Section 4]{ESBH}, and it turns out that it has
 hypersurface singularities.

 Let us give an ``interpretation'' of these singularities via ${\Phi}_{\rm cont}$:
 let $X$ be a classical and exceptional Enriques surface over an algebraically 
 closed field of characteristic $2$.
 Then, there exists a family ${\cal X}\to S$ over some
 local Artinian base with special fiber $X$, as well as a small extension
 $S\to\overline{S}$ such that the family cannot be extended 
 over $\overline{S}$.
 After choosing a Cossec--Verra polarization
 $\cal L$ on $X$, this polarization extends uniquely to ${\cal X}$ 
 (since $h^1(\OO_X)=h^2(\OO_X)=0$), and
 $\Phi_{\rm cont}$ yields a family ${\cal X}'\to S$.
 Since $\ModulPolAmple$ is smooth at $\Phi_{\rm cont}(X,{\cal L})$,
 the family ${\cal X}'\to S$ extends to a family
 $\overline{{\cal X}}'\to \overline{S}$.
 By construction, ${\cal X}\to S$ is a simultaneous resolution 
 of singularities of ${\cal X}'\to S$ over $S$.
 By assumption, a simultaneous resolution of singularities
 of $\overline{{\cal X}}'\to\overline{S}$ extending ${\cal X}\to S$
 does not exist over $\overline{S}$. 
 However, it does exist
 after a {\em ramified} extension of $\overline{S}$ by Artin's result
 \cite{Artin Brieskorn}.

 Summing up, the singularities at $X$ can be explained
 via ${ADE}$-curves and obstructions coming from Artin's simultaneous 
 resolution functor. 
 Over the complex numbers, similar phenomena % for surfaces of general type
 have been described in \cite[Section 4]{Burns; Wahl}.
\end{Remark}


\begin{thebibliography}{XXXXX}
  \bibitem[ACGH]{acgh} E.~Arbarello, M.~Cornalba, P.~Griffiths, J.~Harris,
     {\em Geometry of algebraic curves. Vol. I.},
     Grundlehren der Mathematischen Wissenschaften 267, Springer 1985.   
  \bibitem[Ar62]{Artin Some numerical} M.~Artin, {\em Some numerical criteria for 
     contractability of curves on algebraic surfaces}, Amer. J. Math. 84, 485-496 (1962). 
   \bibitem[Ar74a]{Artin versal deformations} M.~Artin, {\em Versal 
     Deformations and Algebraic Stacks}, Invent. Math. 27, 165-189 (1974).
  \bibitem[Ar74b]{Artin Brieskorn} M.~Artin, {\em Algebraic construction of 
     Brieskorn's resolutions}, J. Algebra 29, 330-348 (1974).
  \bibitem[B-P83]{Barth; Peters} W.~ Barth, C.~Peters, {\em Automorphisms of 
     Enriques surfaces}, Invent. Math.  73, 383-411 (1983).
  \bibitem[BHPV]{bhpv} W. Barth, K. Hulek, C. Peters, A. Van de Ven, {\em Compact complex surfaces}, 
     Second edition, Ergebnisse der Mathematik und ihrer Grenzgebiete 4, Springer (2004). 
  \bibitem[B-M76]{bm3} E.~Bombieri, D.~Mumford, {\em Enriques' classification 
    of surfaces in char.$p$, III}, Invent. Math. 35, 197-232 (1976).
  \bibitem[B-W74]{Burns; Wahl} D.~M.~Burns, J.~M.~Wahl, {\em Local contributions 
    to global deformations of surfaces}, Invent. Math. 26, 67-88 (1974).
  \bibitem[Ca04]{Casnati} G.~Casnati, {\em The moduli space of Enriques surfaces 
    with a polarization of degree $4$ is rational}, Geom. Dedicata  106, 185-194 (2004).
  \bibitem[CFHR]{cfhr} F.~Catanese, M.~Franciosi, K.~Hulek, M.~Reid,
   {\em Embeddings of curves and surfaces}, Nagoya Math. J. 154, 185-220 (1999).
  \bibitem[Co83]{Cossec Reye} F.~R.~Cossec, {\em Reye congruences},
   Trans. Amer. Math. Soc. 280, 737-751 (1983).
  \bibitem[Co85]{Cossec Models} F.~R.~Cossec, {\em Projective Models of Enriques
   Surfaces}, Math. Ann. 265, 283-334 (1983).
  \bibitem[C-D85]{Cossec Dolgachev Automorphisms} F.~R.~Cossec, I.~Dolgachev,
   {\em On automorphisms of nodal Enriques surfaces},
   Bull. Amer. Math. Soc. 12, 247-249 (1985).
  \bibitem[C-D89]{Cossec; Dolgachev} F.~R.~Cossec, I.~Dolgachev, 
    {\em Enriques surfaces I}, Progress in Mathematics 76, Birkh\"auser 1989.
  \bibitem[Do84]{Dolgachev Automorphisms} I.~Dolgachev, {\em On automorphisms 
   of Enriques surfaces},  Invent. Math. 76, 163-177 (1984). 
%  \bibitem[D-I87]{Deligne; Illusie} P.~Deligne, L.~Illusie, {\em Rel\'evements 
%    modulo $p^2$ et d\'ecomposition du complexe de de Rham}, 
%    Invent. Math.  89, 247-270 (1987).
  \bibitem[E-H87]{eh} D.~Eisenbud, J.~Harris, {\it On Varieties of Minimal Degree
     (A Centennial Account)}, Algebraic Geometry, Bowdoin 1985, Proc. Symp. Pure
     Math. 46, Part 1, 3-13 (1987).
%  \bibitem[E-SB]{Ekedahl unpublished} T.~Ekedahl, N.~I.~Shepherd-Barron,
%    {\em Moduli and Periods of Simply Connected Enriques Surfaces},
%     unpublished manuscript.
  \bibitem[E-SB04]{Ekedahl} T.~Ekedahl, N.~I.~Shepherd-Barron, {\em On exceptional 
    Enriques surfaces}, arXiv:math/0405510 (2004).
  \bibitem[E-SB-H12]{ESBH} T.~Ekedahl, N.~I.~Shepherd-Barron, J.~M.~E.~Hyland,
    {\em Moduli and periods of simply connected Enriques surfaces},
    arXiv:1210.0342 (2012).
  \bibitem[El78]{Elkik} R.~Elkik, {\em Singularit\'es rationnelles et d\'eformations},
    Invent. Math. 47, 139-147 (1978).    
  \bibitem[En08]{Enriques} F.~Enriques, {\em Un' osservazione relativa alla 
    superficie di bigenere uno}, Bologna Rend. 12, 40-45 (1908).
%  \bibitem[FGA]{fga} A.~Grothendieck, {\em Fondements de la g\'eom\'etrie alg\'ebrique},
%   Extraits du S\'eminaire Bourbaki, 1957--1962.
%  \bibitem[Gr62]{Grothendieck Picard} A.~Grothendieck, {\em Technique de descente 
%    et th\'eor\`emes d'existence en g\'eom\'etrie alg\'ebrique. VI. 
%    Les sch\'emas de Picard: propri\'et\'es g\'en\'erales},
%    S\'eminaire Bourbaki, t. 14, n$^{\rm o}$ 236, 1961/62.
  \bibitem[Il79]{Illusie} L.~Illusie, {\em Complexe de de Rham--Witt et 
     cohomologie cristalline}, Ann. Sci. \'Ecole Norm. Sup. 12, 501-661 (1979).
  \bibitem[Jo83]{Jouanolou} J.-P. Jouanolou, {\em Th\'eor\`emes de Bertini et 
     applications}, Progress in Mathematics 42, Birkh\"auser 1983.
  \bibitem[Kn71]{Knutson} D.~Knutson, {\em Algebraic spaces}, Lecture Notes in Mathematics 203,
    Springer 1971.  
  \bibitem[Ko94]{Kondo} S.~Kond\={o}, {\em The rationality of the moduli space 
     of Enriques surfaces},  Compositio Math.  91, 159-173 (1994).
  \bibitem[La83]{Lang} W.~E.~Lang, {\em On Enriques surfaces in characteristic 
     $p$. I},  Math. Ann.  265, 45-65 (1983).
%  \bibitem[La88]{Lang2} W.~E.~Lang, {\em On Enriques surfaces in characteristic 
%     $p$. II},  Math. Ann.  281, 671-685 (1988).
  \bibitem[L-M00]{Laumon} G.~Laumon, L. Moret-Bailly, {\em Champs alg\'ebriques}, 
    Ergebnisse der Mathematik und ihrer Grenzgebiete 39, Springer (2000).
  \bibitem[Li08]{Liedtke Horikawa} C.~Liedtke, {\em Algebraic Surfaces with Small $c_1^2$
      in Positive Characteristic}, Nagoya Math. J. 191, 111-134 (2008).
%  \bibitem[Li09]{Liedtke} C.~Liedtke, {\em Non-classical Godeaux Surfaces},
%    Math. Ann. 343, 623-637 (2009).
%  \bibitem[Ma89]{Matsumura} H.~Matsumura, {\em Commutative ring theory},
%   Second edition, Cambridge University Press (1989).
%  \bibitem[M-M64]{Matsusaka} T.~Matsusaka, D.~Mumford, {\em Two fundamental theorems
%      on deformations of polarized varieties}, Amer. J. Math. 86, 668-684 (1964).
%  \bibitem[Ma68]{Matsusaka} T.~Matsusaka, {\em Algebraic deformations of polarized 
%    varieties},  Nagoya Math. J. 31, 185-245 (1968).
 \bibitem[Mu66]{Mumford} D.~Mumford, {\em Lectures on curves on an algebraic surface},
    Annals of Mathematics Studies 59, Princeton University Press 1966.
  \bibitem[O-T70]{Tate; Oort} F.~Oort, J.~Tate, {\em Group schemes of prime order},
    Ann. Sci. \'Ecole Norm. Sup. 3, 1-21 (1970).
  \bibitem[Ra70]{Raynaud} M.~Raynaud, {\em Sp\'ecialisation du foncteur de Picard},
    Inst. Hautes \'Etudes Sci. Publ. Math. 38, 27-76 (1970).
  \bibitem[Ra79]{Raynaud p-torsion} M.~Raynaud, {\em ''$p$-torsion`` du schema de Picard}, 
    Ast\'erisque 64, Soc. Math. de France, 195-280 (1979).    
  \bibitem[Ri96]{rizov} J.~Rizov, {\em Moduli stacks of polarized $K3$ surfaces 
    in mixed characteristic}, Serdica Math. J. 32, 131-178 (2006).
  \bibitem[SD74]{Saint-Donat} B.~Saint-Donat, {\em Projective models of $K3$ 
    surfaces}, Amer. J. Math. 96, 602-639 (1974).
  \bibitem[Sa03]{Salomonsson} P.~Salomonsson, {\em Equations for some very special
    Enriques surfaces in characteristic two}, arxiv:math.AG/0309210 (2003).
%  \bibitem[Sch68]{Schlessinger} M.~Schlessinger, {\em Functors of Artin rings},
%   Trans. Amer. Math. Soc. 130, 208-222 (1968).
  \bibitem[Sch91]{Schreyer} F.-O. Schreyer, {\em A standard basis approach to syzygies 
    of canonical curves}, J. Reine Angew. Math. 421, 83-123 (1991).
  \bibitem[Ve83]{Verra} A.~Verra, {\em The \'etale double covering of an 
     Enriques' surface}, Rend. Sem. Mat. Univ. Politec. Torino 41, 131-167 (1983).
%  \bibitem[Ve83]{Verra} A.~Verra, {\em On Enriques surface as a fourfold cover 
%   of $\PP^2$},  Math. Ann. 266, 241-250 (1983).
\end{thebibliography}
\end{document}